\documentclass[11pt,amsfonts]{article}
\usepackage{graphicx}
\usepackage{latexsym}
\usepackage{amssymb}
\usepackage{amsmath}
\usepackage{enumerate}

\usepackage{amsmath}
\usepackage{stmaryrd}
\usepackage{hyperref}
\usepackage{amssymb}
\usepackage{CJK}
\usepackage{color}

\usepackage{amsmath,latexsym,amssymb,amsthm,array,amsfonts}

\usepackage{mathrsfs,dsfont,bbm}

\usepackage{layout}
\usepackage{eufrak}

\newtheorem{lemma}{Lemma}

\newtheorem{theorem}{Theorem}
\newtheorem{proposition}{Proposition}
\newtheorem{remark}{Remark}
\newtheorem{hypothesis}{Hypothesis}
\newtheorem{assumption}{Assumption}
\newtheorem{example}{Example}

\def\real{{\mathord{{\rm I\kern-2.8pt R}}}}        
\def\inte{{\mathord{{\rm I\kern-2.8pt N}}}}

\def\sZZ{{\rm Z\kern-2.8ptem{}Z}}

\def\z{{\mathchoice
  {\sZZ}
  {\sZZ}
  {\rm Z\kern-0.30em{}Z}
  {\rm Z\kern-0.25em{}Z} }}
\def\sQQ{{\kern 0.27em \vrule height1.45ex width0.03em depth0em
          \kern-0.30em \rm Q}}
\def\qu{{\mathchoice
    {\sQQ}
    {\sQQ}
  {\kern 0.225em \vrule height1.05ex width0.025em depth0em \kern-0.25em \rm Q}
  {\kern 0.180em \vrule height0.78ex width0.020em depth0em \kern-0.20em \rm Q}
        }}
\def\sCC{{\kern 0.27em \vrule height1.45ex width0.03em depth0em
          \kern-0.30em \rm C}}
\def\complex{{\mathchoice
    {\sCC}
    {\sCC}
  {\kern 0.225em \vrule height1.05ex width0.025em depth0em \kern-0.25em \rm C}
  {\kern 0.180em \vrule height0.78ex width0.020em depth0em \kern-0.20em \rm C}
        }}

\newcommand{\ba}{\begin{array}}
\newcommand{\ea}{\end{array}}
\newcommand{\be}{\begin{equation}}
\newcommand{\ee}{\end{equation}}
\newcommand{\bea}{\begin{eqnarray}}
\newcommand{\eea}{\end{eqnarray}}
\newcommand{\beaa}{\begin{eqnarray*}}
\newcommand{\eeaa}{\end{eqnarray*}}

\def\z{\zeta}

\font\tenmath=msbm10 \font\sevenmath=msbm7 \font\fivemath=msbm5
\newfam\mathfam \textfont\mathfam=\tenmath
\scriptfont\mathfam=\sevenmath \scriptscriptfont\mathfam=\fivemath

\def \={{\buildrel {\rm (law)} \over =}}

\def\qed{ \hfill \vrule width.25cm height.25cm depth0cm\smallskip}

\newcommand{\basa}{\begin{assumption}}
\newcommand{\easa}{\end{assumption}}

\newcommand{\bas}{\begin{assum}}
\newcommand{\eas}{\end{assum}}

\newcommand{\ignore}[1]{}
\begin{document}

\renewcommand{\thefootnote}{\fnsymbol{footnote}}

\renewcommand{\thefootnote}{\fnsymbol{footnote}}

\title{Generalized space-time fractional stochastic kinetic equation}

\author{Junfeng Liu\\
School of Statistics and Data Science, \\
Nanjing Audit University, Nanjing, P.R. China.\\
junfengliu@nau.edu.cn \vspace*{0.1in} \\
 }
\maketitle

\begin{abstract}

In this paper, we study a class of nonlinear space-time fractional stochastic kinetic equations in $\mathbb{R}^d$ with Gaussian noise which is white in time and homogeneous in space.  This type of equation  constitutes an extension of the non-linear stochastic heat equation involving fractional derivative in time and fractional Laplacian in space. We give a necessary condition on the spatial covariance  for the existence and uniqueness of the solution. We also study various properties of the solution: path regularity, the behavior of second moment  and the stationarity in the case of linear additive noise.

\end{abstract}

\vskip0.3cm

{\bf 2010 AMS Classification Numbers: } 60G22; 60H07; 60H15

 \vskip0.3cm


\section{Introduction}
Fractional stochastic partial differential equations constitute a subclass of stochastic partial differential equations. The main characteristic of this class of stochastic equations is that they involve  fractional derivatives and integrals, which replace the usual derivatives and integrals.  The fractional stochastic partial differential equations received a particular  attention in the last decades because they emerge in anomalous diffusion models in physich, among other areas of applications.

The aim of the present article is to study the following space-time fractional stochastic kinetic equations, for any $(t,x)\in\mathbb{R}_+\times\mathbb{R}^d$
\begin{equation}\label{sec1-eq1.1}
\left\{
\begin{aligned}
& \left(\frac{\partial^\beta }{\partial t^\beta}+ \nu({\rm I}-\Delta)^{\gamma/2}(-\Delta)^{\alpha/2}\right)u(t,x)=  I_t^{1-\beta}\left(\lambda\sigma(u(t,x))\dot{W}(t,x)\right),\\
&u(0,x)=u_0(x),\quad x\in\mathbb{R}^d,\\
\end{aligned}
\right.
\end{equation}
where $\beta\in(0,1]$, $\gamma\geq0,\alpha>0$ are some fractional parameters, $\nu$ and $\lambda$ are two positive parameters,  $\lambda$ being called the intensity of the noise. The coefficient  $\sigma(\cdot)$ is a measurable function, and $\dot{W}$ is a Gaussian noise, white in time and correlated in space. Here $\Delta$ is the $d$-dimensional Laplace operator and the operators $(I-\Delta)^{\gamma/2},\gamma\geq0$ and $(-\Delta)^{\alpha/2},\alpha>0$ are interpreted as the inverses of the Bessel and Riesz potentials respectively. We will specify later the required conditions on the function $\sigma(\cdot)$ and the Gaussian noise $\dot{W}$. In Eq. \eqref{sec1-eq1.1}, the time derivative operator $\frac{\partial^\beta}{\partial t^\beta}$ with order $\beta\in(0,1]$ is defined in the {\it Caputo-Djrbashian} sense (for example, Caputo~\cite{Caputo1967}, Anh and Leonenko \cite{AnhJSP2001}):
\begin{equation}
\frac{\partial^\beta}{\partial t^\beta} u(t,x)=\left\{
\begin{aligned}
& \frac{1}{\Gamma(1-\beta)}\left[\frac{\partial}{\partial t}\int_0^t\frac{u(s,x)}{(t-s)^{\beta}}ds-\frac{u(0,x)}{t^\beta}\right],\quad {\rm if}\quad \beta\in(0,1),\\
&\frac{\partial}{\partial t}u(t,x),\quad {\rm if} \quad \beta=1.\\
\end{aligned}
\right.
\end{equation}
The fractional integral operator $I_t^{1-\beta}$ is defined by
$$
I_t^{1-\beta} u(t,x):=\frac{1}{\Gamma(1-\beta)}\int_0^t\frac{u(s,x)}{(t-r)^{\beta}}ds,\quad \beta\in(0,1).
$$

The deterministic counterparts of Eq. \eqref{sec1-eq1.1} have received a lot of attention, see for example, Anh and Leonenko \cite{AnhJSP2001}, Caputo \cite{Caputo1967}, Meerschaert {\it et al} \cite{MeerschaertAOP2009}, Nane \cite{Nane2012} and references therein.
This is because they appear to be very useful for modeling,  being  introduced to describe physical phenomena such as diffusion in porous media with fractal geometry, kinematics in viscoelastic media, relaxation processes in complex systems (including viscoelastic materials, glassy materials, synthetic polymers, biopolymers), propagation of seismic waves, anomalous diffusion and turbulence (see Anh and Leonenko \cite{AnhJSP2001}, Caputo \cite{Caputo1967} and references therein).  Such equations are obtained from the classical diffusion equation by replacing the first or second-order derivative by a fractional derivative.


The fractional SPDEs   represent a combination of the deterministic fractional equations and the stochastic integration theory developped by Walsh (see  \cite{Walsh1986}, see also Dalang's seminal paper \cite{DalangEJP1999}). Several types of fractional SPDEs have been considered in  Chen, Kim and Kim \cite{ChenSPA2015},  Chen {\it et al} \cite{ChenStochastics2016}, Foondun and Nane \cite{Foondun2015}, Hu and Hu \cite{Hu2015}, Liu and Yan \cite{LiuJOTP2014}, M\'{a}rquez-Carreras \cite{Carreras2009}, \cite{Carreras2012}, Mijena and Nane \cite{MijenaSPA2015}, \cite{MijenaPA2015} and and  references therein.

 In this work, we are interested in space-time fractional SPDE \eqref{sec1-eq1.1}. It includes some  widely studied particular cases. We refer,  for example, to the classical stochastic heat equation with $\beta=1$, $\gamma=0$ and $\alpha=2$ (see e.g. Dalang \cite{DalangEJP1999}, Khoshnevisan \cite{Khoshnevisan2014}), the fractional stochastic heat equation with $\beta=1$, $\gamma=0$ and $\alpha>0$ (see example, Chen and Dalang \cite{ChenSPDE-AC2015}, \cite{ChenAIHP2015},  Foondun and Nane \cite{Foondun2015}, M\'{a}rquez-Carreras \cite{Carreras2012}, Tudor\cite{Tudor2013}), the generalized fractional kinetic equation with $\beta=0$, $\gamma\geq0$ and $\alpha>0$ (see example M\'{a}rquez-Carreras \cite{Carreras2009}), the space-time fractional stochastic partial differential equation with $0<\beta<1$, $\gamma=0$ and $0<\alpha\leq2$ (see example, Mijena and Nane \cite{MijenaSPA2015}, \cite{MijenaPA2015}).

To be more precise, we extend the result in \cite{MijenaSPA2015}, \cite{Carreras2009} and \cite{FMN} by including in the model the Bessel operator $(I-\Delta) ^ {\frac{\gamma }{2}}$ with $\gamma \geq 0$ and by generalizing the stochastic noise, in the sense that we  allow a more general structure for the spatial covariance of the Gaussian noise $W$ in \eqref{sec1-eq1.1} (which is taken to be space-time white noise in \cite{MijenaSPA2015} and colored by a Riesz kernel in space  in \cite{FMN}).  The presence of this Bessel operator  brings more flexibility to the model, by including for $\gamma=0$ the situation treated in \cite{MijenaSPA2015}, \cite{FMN} or \cite{Carreras2009}. From the technical point of view, the appearance of the Bessel operator leads to a new expression of the fundamental solution associated to the equation  \eqref{sec1-eq1.1}. Indeed, we need new technical estimates for this kernel, which are obtained in Section 2.2. The Bessel operator is also essential in order to get an asymptotically stationary solution, as discussed in Section 4 of our work.  Concretely,   we study the existence and uniqueness of the solution to Eq. \eqref{sec1-eq1.1} under global Lipschitz conditions on diffusion coefficient $\sigma$ by using the random field approach of Walsh \cite{Walsh1986} and time fractional Duhamel's principle (see e.g. Umarov \cite{UmarovJDE2012} and Mijena and Nane \cite{MijenaSPA2015}). Moreover we study some new properties for the solution to time-space fractional SPDE \eqref{sec1-eq1.1},  including an  upper bound of the second moment, the H\"{o}lder regularity in time and space variables, the (asymptotically) stationarity of the solution with respect to time and space variables in some particular case.

We organize this paper as follows. In Section \ref{sec2} we introduce the Gaussian noise $\dot{W}(t,x)$, and we prove some properties of Green function $G_t(x)$ associated with the fractional heat type equation \eqref{Greenfunction}. In Section \ref{sec3}  we give our  main result about existence and uniqueness of the solution and some properties of the solution, including  the H\"{o}lder regularity and the behavior of the second moment.  In Section \ref{sec4} we study the linear additive case, with zero initial condition, i.e.   $u_0(x)\equiv0$ and $\sigma(x)\equiv1$. We  see  that  the solution of \eqref{sec1-eq1.1} is a Gaussian field with zero mean, with  stationary increments, and a continuous covariance function in space, while the it is not stationary in time but tends to a stationary process when the time goes to infinity.

\section{Preliminaries}\label{sec2}

In this section, we recall some basic properties of  the stochastic integral with respect to the Gaussian noise $\dot{W}$ appearing in Eq. \eqref{sec1-eq1.1} and  some basic facts  on the solution to the fractional heat  equation \eqref{Greenfunction}.

\subsection{Gaussian noise}\label{sec2.1}

We denote by $C_0^\infty(\mathbb{R}_+\times\mathbb{R}^{d})$ the space of infinitely differentiable functions on $\mathbb{R}_+\times\mathbb{R}^d$ with compact support and by
$\mathcal{S}(\mathbb{R}^d)$ the Schwartz space of rapidly decreasing $C^\infty$ functions in $\mathbb{R}^d$ and let $\mathcal{S}'(\mathbb{R}^d)$ denote its dual space of rapidly decreasing infinitely differentiable functions on $\mathbb{R}^d$. For $\varphi\in L^1(\mathbb{R}^d)$, we let $\mathcal{F}\varphi$
be the Fourier transform of $\varphi$ defined by
\begin{equation}\label{def-fourier}
\mathcal{F}\varphi(\xi)=\int_{\mathbb{R}^d}e^{-{\rm i}\xi\cdot x}\varphi(x)dx,\quad \xi\in\mathbb{R}^d.
\end{equation}

We begin by introducing the framework in  \cite{DalangEJP1999}. Let $\mu$ be a non-negative tempered measure on $\mathbb{R}^d$, i.e., a non-negative measure which satisfies:
\begin{equation}\label{tempered-measure}
\int_{\mathbb{R}^d}
\left(\frac{1}{1+|\xi|^2}\right)^m\mu(d\xi)<\infty,
\end{equation}
for some $m>0$. Since the integrand is non-increasing in $m$, we may assume that $m\geq1$ is an integer. Note that $1+|\xi|^2$ behaves like a constant around $0$, and like $|\xi|^2$ at $\infty$, and hence \eqref{tempered-measure} is equivalent to
$$
\int_{|\xi|\leq1}\mu(d\xi)<\infty \quad {\rm and}\quad \int_{|\xi|\geq1}\frac{1}{|\xi|^{2m}}\mu(d\xi)<\infty,
$$
for some integer $m\geq1$.

Let $f:\mathbb{R}^d\rightarrow\mathbb{R}_+$ be the Fourier transform of a non-negative tempered measure $\mu$ in $\mathcal{S}'(\mathbb{R}^d)$, that is
$$
\int_{\mathbb{R}^d}f(x)\varphi(x)dx=\int_{\mathbb{R}^d}\mathcal{F}\varphi(\xi)\mu(d\xi), \quad\forall\varphi\in\mathcal{S}(\mathbb{R}^d),
$$
where $\mathcal{F}$ denotes the Fourier transform given by \eqref{def-fourier}.
Simple properties of the Fourier transform yield that for any $\phi, \varphi\in\mathcal{S}(\mathbb{R}^d)$
\begin{equation}\label{cov-1}
\int_{\mathbb{R}^d}\int_{\mathbb{R}^d}\varphi(x)f(x-y)\phi(y)dxdy
=\int_{\mathbb{R}^d}\mathcal{F}\varphi(\xi)\overline{\mathcal{F}\phi(\xi)}\mu(d\xi), \quad \forall\varphi, \phi\in\mathcal{S}(\mathbb{R}^d).
\end{equation}
An approximation argument shows that the previous equality also holds for indicator functions
$\varphi=1_{A}$ and $\phi=1_{B}$ with $A,B\in\mathcal{B}_b(\mathbb{R}^d)$, where $\mathcal{B}_b(\mathbb{R}^d)$ denotes the class of bounded Borel sets of $\mathbb{R}^d$, that is
\begin{equation}\label{fourier-1}
\int_A\int_Bf(x-y)dxdy=\int_{\mathbb{R}^d}\mathcal{F}1_{A}(\xi)\overline{\mathcal{F}1_{B}(\xi)}\mu(d\xi).
\end{equation}

In this article we consider a zero-mean Gaussian process $W=\{W(t,A);t\in[0,T], A\in\mathcal{B}_b(\mathbb{R}^d)\}$ with covariance
$$
\mathbf{E}(W(t,A)W(s,B))=(t\wedge s)\int_A\int_Bf(x-y)dxdy,
$$
on a complete probability space $(\Omega,\mathcal{F}, P)$.

Let $\mathcal{E}$ be the set of linear combinations of elementary functions $\{1_{[0,t]\times A},t\geq0, A\in\mathcal{B}_b(\mathbb{R}^d)\}$. With the Gaussian process $W$, we can associate a canoncial Hilbert space $\mathcal{H}$ which is defined as the closure of $\mathcal{E}$ with respect to the inner product $\langle\cdot,\cdot\rangle_{\mathcal{H}}$ defined by
$$
\langle \varphi, \phi\rangle_{\mathcal{H}}=\int_{\mathbb{R}_+}\int_{\mathbb{R}^d}\int_{\mathbb{R}^d}\varphi(t,x)f(x-y)\phi(t,y)dxdydt.
$$
(Alternatively, $\mathcal{H}$ can be defined as the completion of $C_0^\infty(\mathbb{R}_+\times\mathbb{R}^{d})$ with respect to the inner product $\langle\cdot,\cdot\rangle_{\mathcal{H}}$).

We denote by $W(\varphi)$ the random field indexed by functions $\varphi\in L^2(\mathbb{R}_+\times\mathbb{R}^d)$ and for all $\varphi,\phi\in L^2(\mathbb{R}_+\times\mathbb{R}^d)$, we have
\begin{equation}\label{cov-2}
\begin{split}
\mathbf{E}(W(\varphi)W(\phi))&=\int_{\mathbb{R}_+}\int_{\mathbb{R}^d}\int_{\mathbb{R}^d}
\varphi(t,x)f(x-y)\phi(t,y)dxdydt\\
&=\int_{\mathbb{R}_+}\int_{\mathbb{R}^d}\mathcal{F}\varphi(t,\cdot)(\xi)
\overline{\mathcal{F}\phi(t,\cdot)(\xi)}\mu(d\xi)dt,
\end{split}
\end{equation}
where $\mathcal{F}\varphi(t,\cdot)(\xi)$ denotes the Fourier transform with respect to the space variable of $\varphi(t,x)$ only. Hence $W(\varphi)$ can be represented as
$$
W(\varphi)=\int_{\mathbb{R}_+}\int_{\mathbb{R}^d}\varphi(t,x)W(dx,dt).
$$
Note that $W(\varphi)$ is $\mathcal{F}_t$-measurable whenever $\varphi$ is supported on $[0,t]\times\mathbb{R}^d$.

\begin{remark}\label{remark1}
Since the spectral measure $\mu$ is non-trivial positive tempered measure,
we can ensure that there exist positive constants $c_1,c_2$ and $k$ such that
\begin{equation}\label{est-mu}
c_1<\int_{\{|\xi|<k\}}\mu(d\xi)<c_2.
\end{equation}
\end{remark}

As usual, the Gaussian process $W$ can be extended to a {\it worthy martingale measure}, in the sense given by Walsh \cite{Walsh1986}. Dalang \cite{DalangEJP1999} presented an extension of Walsh's stochastic integral that requires the following integrability condition in terms of the Fourier transform of $G$
\begin{equation}\label{DC}
\int_0^Tdt\int_{\mathbb{R}^d}\mu(d\xi)|\mathcal{F}G_t(\cdot)(\xi)|^2<\infty,
\end{equation}
where $G$ is the fundamental solution of
\begin{equation} \label{Greenfunction}
\left(\frac{\partial^\beta }{\partial t^\beta}+ \nu(I-\Delta)^{\gamma/2}(-\Delta)^{\alpha/2}\right)G_t(x)=0.
\end{equation}
Provided that \eqref{DC} is satisfied and assuming conditions on $\sigma(\cdot)$ that will be described later, following Walsh \cite{Walsh1986}, we will understand a solution of \eqref{sec1-eq1.1} to be a jointly measurable adapted process $\{u(t,x), (t,x)\in[0,T]\times\mathbb{R}^d\}$ satisfying the integral equation
\begin{equation}\label{mildsolution}
u(t,x)=(\mathcal{G}u_0)_t(x)+\lambda\int_0^t\int_{\mathbb{R}^d}
G_{t-s}(x-y)\sigma(u(s,y))W(ds,dy),
\end{equation}
where
$$
(\mathcal{G}u_0)_t(x)=\int_{\mathbb{R}^d}G_t(x-y)u_0(y)dy,
$$
and the stochastic integral in \eqref{mildsolution} is defined with respect to the $\mathscr{F}$-martingale measure $W(t,A)$. Next we give the meaning  of Walsh-Dalang integrals that is used in \eqref{mildsolution}. (For the details, we refer the readers to Dalang \cite{DalangEJP1999}).

\begin{enumerate}

\item We say that  $(t,x)\rightarrow\Phi_t(x)$ is an elementary random field when there exist $ 0\leq a<b$,  a  $\mathcal{F}_a$-measurable random variable $X\in L^2(\Omega)$ and a deterministioc function $\phi\in L^2(\mathbb{R}^d)$ such that
    $$
    \Phi_t(x)=X1_{[a,b]}(t)\phi(x),\quad t>0,x\in\mathbb{R}^d.
    $$
\item If $h=h_t(x)$ is non-random and $\Phi$ is elementary as above, then we set
\begin{equation}\label{n11}
\int h\Phi dW:=X\int_{(a,b)\times\mathbb{R}^d}h_t(x)\phi(x)W(dt,dx).
\end{equation}
\item The stochastic integral in (\ref{n11}) is a  Wiener integral, and it is well defined if and only if $h_t(x)\phi(x)\in L^2([a,b]\times\mathbb{R}^d)$.
\item Under the above notation, we have the Walsh isometry
$$
\mathbf{E}\left(\left|\int h\Phi dW\right|^2\right)=\int_0^Tds\int_{\mathbb{R}^d}dyh_s(y)^2
\mathbf{E}(|\Phi_s(y)|^2).
$$

\end{enumerate}

\subsection{Some properties of the fundamental solution}
We will give some estimates for the fundamental solution associated  to  (\ref{sec1-eq1.1}). The properties of this fundamental solution will play an important role in the sequel.

Let $G_t(x)$ be the fundamental solution of the fractional kinetic equation \eqref{Greenfunction} with $\beta\in (0,1],\nu>0$, and $\gamma\geq0,\alpha>0$. Anh and Leonenko \cite{AnhJSP2001}, showed that, \eqref{Greenfunction} is equivalent to
the Cauchy problem
\begin{equation}\label{Cauchyproblem}
(\mathscr{D}_t^\beta\mathcal{F}G_t(\cdot))(\xi)
+\nu|\xi|^\alpha(1+|\xi|^2)^{\gamma/2}\mathcal{F}G_t(\cdot)(\xi)=0,\quad
\mathcal{F}G_0(\cdot))(\xi)=1
\end{equation}
abnd they also have proved that \eqref{Cauchyproblem} has a unique solution given by
\begin{equation}\label{Fourier-Greenfunction}
\mathcal{F}G_t(\cdot)(\xi)=E_\beta(-\nu t^\beta|\xi|^\alpha(1+|\xi|^2)^{\gamma/2}),\quad \beta>0,
\end{equation}
where
\begin{equation}\label{Def-E-1}
E_\beta(x) =\sum_{j=0}^\infty\frac{x^j}{\Gamma(1+\beta j)},\quad x>0,
\end{equation}
is the Mittag-Leffler function of order $\beta$. The inverse Fourier transform yields that
\begin{equation}\label{Greenfunction-1}
G_t(x)=(2\pi)^{-d}\int_{\mathbb{R}^d}e^{{\rm i}\langle \xi,x\rangle}E_\beta(-\nu t^\beta|\xi|^\alpha(1+|\xi|^2)^{\gamma/2})d\xi.
\end{equation}

We know that
\begin{equation}\label{e-beta in L1}
E_\beta(-\nu t^\beta|\xi|^\alpha(1+|\xi|^2)^{\gamma/2})\in L^1(\mathbb{R}^d),
\end{equation}
for every $0<\beta\leq1$ if $\alpha+\gamma>d$. From this range we see the role played of the parameter $\gamma$ in Eq. \eqref{Fourier-Greenfunction}.

Moreover one has the uniform estimates of Mittag-Leffler function (e.g. Theorem 4 in Simon \cite{SimonEJP2014})
\begin{equation}\label{estimates-E-beta}
\frac{1}{1+\Gamma(1-\beta)x}\leq E_\beta(-x)\leq\frac{1}{1+\Gamma(1+\beta)^{-1}x}, \quad {\rm for}\quad x>0.
\end{equation}

The following lemma gives a sharp estimate for the  $L^ {2}$ norm (in time) of the Green kernel. It extends Lemma 1 in \cite{MijenaSPA2015} and Lemma 2.1 in \cite{FMN}.

\begin{lemma}\label{est-G-square}
For $0<\beta<1$ and $d<2(\alpha+\gamma)$, we have the following
\begin{equation}\label{g-square}
\int_{\mathbb{R}^d}G_t(x)^2dx\leq  C_2t^{-\frac{\beta d}{\alpha+\gamma}}.
\end{equation}
where $B\left(\frac{d}{\alpha+\gamma},2-\frac{d}{\alpha+\gamma}\right)$ is a Beta function. The  (strictly) positive constant  is  given by $C_2=\frac{B\left(\frac{d}{\alpha+\gamma},2-\frac{d}{\alpha+\gamma}\right)}{\alpha+\gamma}
\left(\frac{\Gamma(1-\beta)}{\nu}\right)^{\frac{d}{\alpha+\gamma}}
\frac{2\pi^{d/2}}{\Gamma(d/2)}\frac{1}{(2\pi)^{d}}$.
\end{lemma}
\begin{proof}
Using the Plancherel identity  and the equality \eqref{Fourier-Greenfunction}, we can write
\begin{align*}
\int_{\mathbb{R}^d}G_t(x)^2dx&=\frac{1}{(2\pi)^d}\int_{\mathbb{R}^d}|
\mathcal{F}G_t(\cdot)(\xi)|^2d\xi\\
&=\frac{1}{(2\pi)^d}\int_{\mathbb{R}^d}
\left|E_\beta\left(-\nu t^\beta|\xi|^\alpha(1+|\xi|^2)^{\gamma/2}\right)\right|^2d\xi\\
&=\frac{2\pi^{d/2}}{\Gamma(d/2)}\frac{1}{(2\pi)^{d}}\int_{0}^{+\infty}r^{d-1}
\left(E_\beta\left(-\nu t^\beta r^\alpha(1+r^2)^{\gamma/2}\right)\right)^2dr,
\end{align*}
where we have used the integration in polar coordinates in the last equation above and the positive constant resulting from the integration over the angular spherical coordinates. Now using the upper bound in \eqref{estimates-E-beta} and the fact
$r^\alpha(1+r^2)^{\gamma/2}\geq r^{\alpha+\gamma}$ with $r>0$, we get with the change of variable formula $z=\Gamma(1+\beta)^{-1}\nu t^\beta r^{\alpha+\gamma}$,
\begin{align*}
\int_{0}^{+\infty}r^{d-1}
\left(E_\beta\left(-\nu t^\beta r^\alpha(1+r^2)^{\gamma/2}\right)\right)^2dr
&\leq\int_{0}^{+\infty}r^{d-1}\frac{1}{(1+\Gamma(1+\beta)^{-1}\nu t^\beta r^{\alpha+\gamma})^2}dr\\
&=\frac{1}{\alpha+\gamma}
\left(\frac{\Gamma(1+\beta)}{\nu}\right)^{\frac{d}{\alpha+\gamma}}t^{-\frac{\beta d}{\alpha+\gamma}}\int_0^\infty z^{\frac{d}{\alpha+\gamma}-1}(1+z)^{-2}dz.
\end{align*}
Hence $\int_0^\infty z^{\frac{d}{\alpha+\gamma}}(1+z)^{-2}dz<\infty$ if and only if $d<2(\alpha+\gamma)$. In this case, we have
$$
\int_0^\infty z^{\frac{d}{\alpha+\gamma}-1}(1+z)^{-2}dz
=B\left(\frac{d}{\alpha+\gamma},2-\frac{d}{\alpha+\gamma}\right),
$$
where $B\left(\frac{d}{\alpha+\gamma},2-\frac{d}{\alpha+\gamma}\right)$ is a Beta function. Then we can conclude the proof of upper bound in \eqref{g-square}.
\end{proof}

\qed

We now prove \eqref{DC} under an integrability condition on the spectral measure $\mu$ given as follows.

\begin{hypothesis}\label{hypothesis-1}
 Assume that the spectral measure $\mu$ associated to the Gaussian noise $\dot{W}$ satisfies
\begin{equation}\label{DC-extension}
\int_{\mathbb{R}^d}\left(\frac{1}{1+|\xi|^2}\right)^{\varrho}\mu(d\xi)<\infty,
\end{equation}
with the parameter $\varrho$ satisfying
\begin{equation}\label{hypo-1}
\varrho=\left\{
\begin{aligned}
&\alpha+\gamma, \quad {\rm if}\quad  0<\beta<\frac12,\\
&\frac{\alpha+\gamma}{2},\quad {\rm if} \quad \beta=\frac12,\\
&\frac{\alpha+\gamma}{2\beta}, \quad {\rm if}\quad  \frac12<\beta<1.\\
\end{aligned}
\right.
\end{equation}
\end{hypothesis}

\begin{remark}
If the parameter $\beta=1$, Eq. \eqref{sec1-eq1.1} reduces to the SPDE (1.1) studied by M\'{a}rquez-Carreras \cite{Carreras2009}. In this paper it is assumed (\ref{DC-extension}) with  $\varrho=\frac{\alpha+\gamma}{2}$. So when $\beta$ is close to one, the exponent $\varrho$ in \eqref{hypo-1} is coincides with the exponent studied in Lemma 2.1 in M\'{a}rquez-Carreras \cite{Carreras2009}. On the other hand, our assumption (\ref{DC-extension}) is weaker when $\beta$ is close to zero.

\end{remark}

Let us now recall some of the main examples of spatial covariances for the noise  which will be our guiding examples in the remainder of the present paper. Below we denote by $|x|$ the Euclidean norm of $x\in\mathbb{R}^d$.

\begin{example}
{\rm
Let $f(x)=\prod\limits_{i=1}^d H_i(2H_i-1)|x_i|^{2H_i-2}$ with $1/2<H_i<1$ for $i=1,\ldots,d$. Then
$
\mu(d\xi)=\prod_{i=1}^dH_i(2H_i-1)|\xi_i|^{-2H_i+1}d\xi.
$
So \eqref{DC-extension} is equivalent to $\sum\limits_{i=1}^d(2H_i-1)>d-2(\alpha+\gamma)$ if $0<\beta<1/2$, it  is equivalent to $\sum\limits_{i=1}^d(2H_i-1)>d-(\alpha+\gamma)$ if $\beta=1/2$ and when $1/2<\beta<1$, the condition  \eqref{DC-extension} is equivalent to $\sum\limits_{i=1}^d(2H_i-1)>d-\frac{\alpha+\gamma}{\beta}$.
}
\end{example}

\begin{example}\label{example-riesz}
{\rm
Let $f(x)=\gamma_{\delta,d}=|x|^{-(d-\delta)}$ be the Riesz kernel of order $\delta \in(0,d)$, then $\mu(d\xi)=|\xi|^{-\delta}d\xi$ and \eqref{DC-extension} is equivalent to $2(\alpha+\gamma)+\delta>d$ if $0<\beta<1/2$, \eqref{DC-extension} is equivalent to $(\alpha+\gamma)+\delta>d$ if $\beta=1/2$ and \eqref{DC-extension} is equivalent to $\frac{\alpha+\gamma}{\beta}+\delta>d$ if $1/2<\beta<1$.
}
This example is also considered in \cite{FMN}. Their condition (used in Theorem 1.3 in this reference) reads  $\frac{\alpha}{\beta} + \delta >d$. Our assumption (\ref{DC-extension}) gives more flexibility when $\beta$ is close to zero but as well as for $\beta $ close to 1 (because of the new parameter $\gamma$ in the expression of the Bessel operator $(I-\Delta)  ^ {\frac{\gamma}{2}}.$
\end{example}

\begin{example}\label{example-bessel}
{\rm
For the Bessel kernel of order $\tau>0$ given by
$
f(x)=\gamma_\tau\int_0^\infty\omega^{\frac{\tau-d}{2}-1} e^{-\omega}e^{-\frac{|x|^2}{4\omega}}d\omega.
$
Then $\mu(d\xi)=(1+|\xi|^2)^{-\frac{\tau}{2}}d\xi$. So \eqref{DC-extension} is equivalent to $2(\alpha+\gamma)+\tau>d$ if $0<\beta< 1/2$, \eqref{DC-extension} is equivalent to $(\alpha+\gamma)+\tau>d$ if $\beta=1/2$ and condition \eqref{DC-extension} is equivalent to $\frac{\alpha +\gamma}{\beta}+\tau>d$ if $1/2<\beta<1$.
}
\end{example}

\begin{example}
{\rm

Let $f(0)<\infty$ (i.e. $\mu$ is a finite measure). It corresponds to a spatially smooth noise $\dot{W}$.
}
\end{example}

\begin{example}
{\rm

Suppose $d=1$ and $f=\delta_0$ (i.e. $\mu$ is the Lebesgue measure). This corresponds to a (rougher) noise $\dot{W}$ which is white in spatial variable.
}
\end{example}

 For any $t\in\mathbb{R}_+$, denote by
\begin{equation}\label{nt-xi}
N_t(\xi)=\int_0^t|\mathcal{F}G_u(\cdot)(\xi)|^2du.
\end{equation}
Then we have the following
\begin{proposition}\label{nt-xi-lower and upper}
Assuming that $t\in\mathbb{R}_+$ and $\xi\in\mathbb{R}^d$, there exist (strictly) positive constants $C_{2.i}(t),i=1,2,3,4$ (depending on $t$) such that
\begin{equation}\label{ests-nt-xi-1}
N_t(\xi)\leq C_{2.2}(t)\left(\frac{1}{1+|\xi|^2}\right)
^{\alpha+\gamma}, \quad {\rm if}\quad 0<\beta<1/2,
\end{equation}

\begin{equation}\label{ests-nt-xi-2}
 N_t(\xi)\leq C_{2.3}(t)\left(\frac{1}{1+|\xi|^2}\right)
^{\frac{\alpha+\gamma}{2}}, \quad {\rm if}\quad \beta=1/2,
\end{equation}
and
\begin{equation}\label{ests-nt-xi-3}
 N_t(\xi)\leq C_{2.4}(t)\left(\frac{1}{1+|\xi|^2}\right)
^{\frac{\alpha+\gamma}{2\beta}}, \quad {\rm if}\quad 1/2<\beta<1.
\end{equation}
The "constants" are defined as follows
\begin{align*}
C_{2.2}(t)&=t+\frac{2^{\alpha+\gamma}\Gamma(1+\beta)^2}{\nu^2(1-2\beta)}t^{1-2\beta},\\
C_{2.3}(t)&=t+2\nu^{-1}\Gamma(3/2)2^{\frac{\alpha+\gamma}{2}}t^{1/2},\\
C_{2.4}(t)&=t+\frac{1}{2\beta-1}\Gamma(1+\beta)^{1/\beta}
\nu^{-1/\beta}2^{\frac{\alpha+\gamma}{2\beta}}.\\
\end{align*}

\end{proposition}

\begin{proof}
For any $t\in\mathbb{R}_+$, from Eqs. \eqref{Fourier-Greenfunction} and \eqref{nt-xi}, we can rewrite $N_t(\xi)$ defined by \eqref{nt-xi} as
$$
N_t(\xi)=\int_0^t\left|E_\beta(-\nu u^\beta|\xi|^\alpha(1+|\xi|^2)^{\gamma/2})\right|^2du
$$
We firstly prove the upper bound for $N_t(\xi)$. By using the upper bound in \eqref{estimates-E-beta} and change of variable $x=\nu u^\beta|\xi|^\alpha(1+|\xi|^2)^{\gamma/2}$, one gets
\begin{align*}
N_t(\xi)&=\frac1\beta\left(\frac{1}{\nu |\xi|^\alpha(1+|\xi|^2)^{\gamma/2}}\right)^{\frac1\beta}\int_0^{\nu t^\beta|\xi|^\alpha(1+|\xi|^2)^{\gamma/2}}
x^{\frac1\beta-1}E_\beta^2(-x)dx.
\end{align*}
We will divide into two cases to estimate it according to the value of $|\xi|$. If $|\xi|\leq 1$, we claim that
\begin{align*}
N_t(\xi)1_{|\xi|\leq 1}&\leq\frac1\beta\left(\frac{1}{\nu |\xi|^\alpha(1+|\xi|^2)^{\gamma/2}}\right)^{\frac1\beta}\int_0^{\nu t^\beta|\xi|^\alpha(1+|\xi|^2)^{\gamma/2}}
x^{\frac1\beta-1}\left(\frac{1}{1+\Gamma(1+\beta)^{-1}x}\right)^2dx\\
&=t.
\end{align*}
If $|\xi|> 1$ and $1/2<\beta<1$ we have
\begin{align*}
N_t(\xi)1_{|\xi|>1}
&\leq\frac1\beta\left(\frac{1}{\nu |\xi|^\alpha(1+|\xi|^2)^{\gamma/2}}\right)^{\frac1\beta}\int_0^{\nu t^\beta|\xi|^\alpha(1+|\xi|^2)^{\gamma/2}}
x^{\frac1\beta-1}\left(\frac{1}{1+\Gamma(1+\beta)^{-1}x}\right)^2dx\\
&=\frac1\beta\left(\frac{\Gamma(1+\beta)}{\nu |\xi|^\alpha(1+|\xi|^2)^{\gamma/2}}\right)^{\frac1\beta}\int_0
^{\nu\Gamma(1+\beta)^{-1} t^\beta|\xi|^\alpha(1+|\xi|^2)^{\gamma/2}}
x^{\frac1\beta-1}\left(1+x\right)^{-2}dx\\
&\leq \frac1\beta\left(\frac{\Gamma(1+\beta)}{\nu |\xi|^\alpha(1+|\xi|^2)^{\gamma/2}}\right)^{\frac1\beta}\int_0
^{\nu\Gamma(1+\beta)^{-1} t^\beta|\xi|^\alpha(1+|\xi|^2)^{\gamma/2}}
\left(1+x\right)^{\frac1\beta-3}dx\\
&\leq\frac{1}{2\beta-1}\left(\frac{\Gamma(1+\beta)}{\nu}\right)
^{\frac1\beta}
\left(\frac{1}{|\xi|^\alpha(1+|\xi|^2)^{\gamma/2}}\right)
^{\frac1\beta}\\
&\leq\frac{1}{2\beta-1}\left(\frac{\Gamma(1+\beta)}{\nu}\right)
^{\frac1\beta}2^{\frac{\alpha+\gamma}{2\beta}}
\left(\frac{1}{1+|\xi|^{2}}\right)^{\frac{\alpha+\gamma}{2\beta}}.
\end{align*}
On the other hand, with $|\xi|>1$ and $0<\beta<1/2$, one gets
\begin{align*}
N_t(\xi)1_{|\xi|>1}
&\leq \frac1\beta\left(\frac{\Gamma(1+\beta)}{\nu |\xi|^\alpha(1+|\xi|^2)^{\gamma/2}}\right)^{\frac1\beta}\int_0
^{\nu\Gamma(1+\beta)^{-1} t^\beta|\xi|^\alpha(1+|\xi|^2)^{\gamma/2}}
x^{\frac1\beta-3}dx\\
&=\frac{\Gamma(1+\beta)^2}{\nu^2(1-2\beta)}t^{1-2\beta}
\left(\frac{1}{\nu|\xi|^\alpha(1+|\xi|^2)^{\gamma/2}}\right)^2\\
&\leq\frac{2^{\alpha+\gamma}\Gamma(1+\beta)^{2}}{\nu^{2}(1-2\beta)}t^{1-2\beta}
\left(\frac{1}{1+|\xi|^{2}}\right)^{\alpha+\gamma}.
\end{align*}
For the critical case $\beta=1/2$, one obtains that
\begin{align*}
N_t(\xi)1_{|\xi|>1}
&\leq 2\left(\frac{\Gamma(3/2)}{\nu |\xi|^\alpha(1+|\xi|^2)^{\gamma/2}}\right)^{2}\int_0
^{\nu\Gamma(3/2)^{-1} t^{1/2}|\xi|^\alpha(1+|\xi|^2)^{\gamma/2}}
x(1+x)^{-2}dx\\
&\leq 2\left(\frac{\Gamma(3/2)}{\nu |\xi|^\alpha(1+|\xi|^2)^{\gamma/2}}\right)^{2}
\int_0
^{\nu\Gamma(3/2)^{-1} t^{1/2}|\xi|^\alpha(1+|\xi|^2)^{\gamma/2}}
(1+x)^{-1}dx\\
&=2\left(\frac{\Gamma(3/2)}{\nu |\xi|^\alpha(1+|\xi|^2)^{\gamma/2}}\right)^{2}\ln\left(1+\nu\Gamma(3/2)^{-1} t^{1/2}|\xi|^\alpha(1+|\xi|^2)^{\gamma/2}\right)\\
&\leq 2\nu^{-1}\Gamma(3/2)t^{1/2}2^{\frac{\alpha+\gamma}{2}}
\left(\frac{1}{1+|\xi|^{2}}\right)^{\frac{\alpha+\gamma}{2}}.
\end{align*}
Then combining the above estimates for $N_t(\xi)$ with $|\xi|\leq1$ and $|\xi|>1$, respectively, we can conclude the proof of bounds  \eqref{ests-nt-xi-1},\eqref{ests-nt-xi-2} and \eqref{ests-nt-xi-3}.
\end{proof}\qed

\begin{remark}

From the above result we see that  Hypothesis \ref{hypothesis-1} implies condition \eqref{DC}. In particular, the estimates  \eqref{ests-nt-xi-1},\eqref{ests-nt-xi-2} and \eqref{ests-nt-xi-3} give the existence of the solution in the linear additive noise cas ($\sigma =1$).
\end{remark}


\section{Existence and uniqueness}\label{sec3}

In this section we will prove the existence and uniqueness of the mild solution to Eq. \eqref{mildsolution}. We first introduce a stronger integrability condition on the spectral measure $\mu$ than Hypothesis \ref{hypothesis-1}. While the existence and uniqueness of the solution can be obtained under Hypothesis 1, the new  assumption presented below will be needed in order to prove certain properties of the solution.

\begin{hypothesis}\label{hypothesis-2}
Assume that the spectral measure $\mu$ associated to $\dot{W}$ satisfies
\begin{equation}\label{DC-extension-2}
\int_{\mathbb{R}^d}\left(\frac{1}{1+|\xi|^2}\right)^{\eta}\mu(d\xi)<\infty,
\end{equation}
with some parameter $\eta$ satisfying
\begin{equation}
\eta\in\left\{
\begin{aligned}
&(0,\alpha+\gamma), \quad {\rm if}\quad  0<\beta<\frac12\\
&(0,\frac{\alpha+\gamma}{2}), \quad {\rm if}\quad  \beta=\frac12,\\
&(0,\frac{\alpha+\gamma }{2\beta}), \quad {\rm if}\quad  \frac12<\beta<1.\\
\end{aligned}
\right.
\end{equation}
\end{hypothesis}

We will need the following estimates for the  Green function  given by \eqref{Greenfunction-1} (their proof is given in the appendix).

\begin{proposition}\label{est-temporal and spatial increments}
Suppose $\beta\in(0,1)$, then we have the following estimates for the temporal and spatial increments of the Green function $G_t(x)$ given by \eqref{Greenfunction-1}
\begin{enumerate}
\item Under Hypothesis \ref{hypothesis-1}, for any $t,t'\in\mathbb{R}_+$ such that $t'<t$ and $x\in\mathbb{R}^d$, we have
\begin{equation}\label{est-time-1}
\int_0^{t'}ds\int_{\mathbb{R}^d}\mu(d\xi)
\left|\mathcal{F}G_{t-s}(x-\cdot)(\xi)-\mathcal{F}G_{t'-s}(x-\cdot)(\xi)\right|^2\leq C_{3.1}|t-t'|^{2\beta},
\end{equation}
with $C_{3.1}=t^{1-2\beta}\int_{|\xi|\leq1}\mu(d\xi)+t^{-2\beta}\int_{|\xi|>1}N_{t'}(\xi)\mu(d\xi)$.

\item Under Hypothesis \ref{hypothesis-2}, for any $t,t'\in\mathbb{R}_+$ such that $t'<t$ and $x\in\mathbb{R}^d$,  we have
\begin{equation}\label{est-time-2}
\int_{t'}^tds\int_{\mathbb{R}^d}\mu(d\xi)
\left|\mathcal{F}G_{t-s}(x-\cdot)(\xi)\right|^2
\leq\left\{
\begin{aligned}
&C_{3.2}|t-t'|^{1-2\beta}, \quad {\rm if}\quad  0<\beta<\frac12,\\
&C_{3.3}|t-t'|^{\frac12}, \quad {\rm if}\quad  \beta=\frac12,\\
&C_{3.4}|t-t'|, \quad {\rm if}\quad  \frac12<\beta<1.\\
\end{aligned}
\right.
\end{equation}
with $C_{3.2}=|t-t'|^{2\beta}\int_{|\xi|\leq1}\mu(d\xi)+\frac{\Gamma(1+\beta)^22^{\alpha+\gamma}}{\nu^2(1-2\beta)}
\int_{|\xi|>1}\left(\frac{1}{1+|\xi|^2}\right)
^{\alpha+\gamma}\mu(d\xi)$, $C_{3.3}=\int_{|\xi|\leq1}\mu(d\xi)|t-t'|^{\frac12}+2^{1+\frac{\alpha+\gamma}{2}}\frac{\Gamma(3/2)}{\nu}\int_{|\xi|>1}\left(\frac{1}{1+|\xi|^2}\right)^{\frac{\alpha+\gamma}{2}}\mu(d\xi)$ \\ and
$C_{3.4}=\int_{|\xi|\leq1}\mu(d\xi)+\frac{c}{2\beta-1}\left(\frac{\Gamma(1+\beta)}{\nu}\right)
^{\frac1\beta}2^{\frac{\alpha+\gamma}{2\beta}}
\int_{|\xi|>1}\left(\frac{1}{1+|\xi|^2}\right)
^{\frac{\alpha+\gamma}{2\beta}}\mu(d\xi)$.
\item Under Hypothesis \ref{hypothesis-2},  for any $t\in\mathbb{R}_+$ and $x,x'\in\mathbb{R}^d$, $\rho_1\in(0,\alpha+\gamma-\eta)$, $\rho_2\in(0,\frac{\alpha+\gamma}{2}-\eta)$ and $\rho_3\in(0,\frac{\alpha+\gamma}{2\beta}-\eta)$  we have
\begin{equation}\label{est-space-1}
\begin{split}
\int_0^tds\int_{\mathbb{R}^d}\mu(d\xi)&
\left|\mathcal{F}G_{t-s}(x-\cdot)(\xi)-\mathcal{F}G_{t-s}(x'-\cdot)(\xi)\right|^2\\
&\leq\left\{
\begin{aligned}
&C_{3.5}|x'-x|^{2\rho_1}, \quad {\rm if}\quad  0<\beta<\frac12,\\
&C_{3.6}|x'-x|^{2\rho_2}, \quad {\rm if}\quad  \beta=\frac12,\\
&C_{3.7}|x'-x|^{2\rho_3}, \quad {\rm if}\quad  \frac12<\beta<1.\\
\end{aligned}
\right.
\end{split}
\end{equation}
with $C_{3.5}=Ct\int_{|\xi|\leq1}\mu(d\xi)+C\frac{t^{1-2\beta}}{1-2\beta}\left( \frac{\Gamma(1+\beta)}{\nu}\right)^2 2^{\alpha+\gamma-\rho_1}
\int_{|\xi|>1}\left(\frac{1}{1+|\xi|^2}\right)
^{\alpha+\gamma-\rho_1}\mu(d\xi)$, $C_{3.6}=Ct\int_{|\xi|\leq1}\mu(d\xi)+C2^{1+\frac{\alpha}{2}-\rho_2}t^{\frac12}\frac{\Gamma(3/2)}{\nu}\int_{|\xi|>1}\left(\frac{1}{1+|\xi|^2}\right)
^{\frac{\alpha+\gamma}{2}-\rho_2}\mu(d\xi)$ and \\
$C_{3.7}=Ct\int_{|\xi|\leq1}\mu(d\xi)+C\frac{1}{2\beta-1}
\left(\frac{\Gamma(1+\beta)}{\nu}\right)^{\frac1\beta}
2^{\frac{\alpha+\gamma}{2\beta}-\rho_3}
\int_{|\xi|>1}\left(\frac{1}{1+|\xi|^2}\right)
^{\frac{\alpha+\gamma}{2\beta}-\rho_3}\mu(d\xi)$.
Notice that all the constants depend on $t$ although we omitt it in the notation.
\end{enumerate}
\end{proposition}

\begin{remark}
\begin{enumerate}
\item Our results of Proposition \ref{est-temporal and spatial increments} extends the results in Mijena and Nane \cite{MijenaSPA2015} to the space-time fractional SPDE with colored Gaussian noises and Khoshnevisan \cite{Khoshnevisan2014} to space-time fractional SPDE, respectively.
\item The above Proposition \ref{est-temporal and spatial increments} here also extends the results in
M\'{a}rquez-Carreras \cite{Carreras2009}, \cite{Carreras2012} to space-time fractional kinetic equation with colored Gaussian noise.
\end{enumerate}

\end{remark}

Let us introduce some additional conditions that we need in order to prove our main results.
The first condition is required for the existence-uniqueness result as well as for
the upper bound on the second moment of the solution.

\begin{assumption}\label{assumption-1}
\begin{enumerate}
\item We assume that the initial condition is a non-random bounded non-negative function $u_0:\mathbb{R}^d\rightarrow\mathbb{R}$.
\item We assume that $\sigma:\mathbb{R}^{d} \rightarrow\mathbb{R} ^{d} $ is Lipschitz continuous satisfying $|\sigma(x)|\leq L_\sigma|x|$ with $L_\sigma$ being a positive constant.
    Moreover for all $x,y\in\mathbb{R}^{d}$
\begin{equation}
|\sigma(x)-\sigma(y)|\leq L_\sigma|x-y|.
\end{equation}
We may assume, with loss of generality, that $L_\sigma$ is also greater than $\sigma(0)$. Since $|\sigma(x)|\leq|\sigma(0)|+L_\sigma|x|$, it follows that $|\sigma(x)|\leq L_\sigma(1+|x|)$ for all $x\in\mathbb{R}^{d}$.
\end{enumerate}
\end{assumption}


Now we can prove the existence and uniqueness of mild solution of Eq.~\eqref{sec1-eq1.1} given by \eqref{mildsolution}.
\begin{theorem}\label{existence}
Under Assumption \ref{assumption-1} and assuming that the spectral measure $\mu$ satisfies Hypothesis \ref{hypothesis-1}, then Eq. \eqref{mildsolution} has a unique adapted solution and for any $t\in\mathbb{R}_+$ and $p\geq1$,
$$
\sup_{(t,x)\in\mathbb{R}_+\times\mathbb{R}^d}\mathbf{E}(|u(t,x)|^p)<\infty.
$$
Moreover, this unique solution is mean-square continuous.
\end{theorem}

\begin{proof}
The proof of existence and uniqueness is standard based on Picard's iterations. For more information, see e.g. Walsh \cite{Walsh1986}, Dalang \cite{DalangEJP1999}. We give a sketch of the proof. Define
\begin{equation}\label{iterated-sequences}
\begin{split}
u^{(0)}(t,x)&=(\mathcal{G}u_0)_t(x),\\
u^{(n+1)}(t,x)&=(\mathcal{G}u_0)_t(x)+\lambda\int_0^t
\int_{\mathbb{R}^d}G_{t-s}(x-y)\sigma(u^{(n)}(s,y))W(ds,dy),\quad n\geq0.
\end{split}
\end{equation}
We could easily prove that the sequence  $\{u^{(n+1)}(t,x),n\geq0\}$ is well-defined and then using Burkholder's inequality, we can show that, for any $n\geq0$ and $t\in\mathbb{R}_+$,
\begin{equation}\label{est-un-lp-1}
\sup_{(t,x)\in\mathbb{R}_+\times\mathbb{R}^d}\mathbf{E}(|u^{(n+1)}(t,x)|^2)<\infty.
\end{equation}
Moreover, by using an extension of Gronwall's lemma (for example, see Lemma 15 in Dalang \cite{DalangEJP1999}), that
\begin{equation}\label{est-un-lp-2}
\sup_{n\geq0}\sup_{(t,x)\in\mathbb{R}_+\times\mathbb{R}^d}\mathbf{E}(|u^{(n+1)}(t,x)|^2)<\infty.
\end{equation}
The same kind of arguments allow us to check \eqref{est-un-lp-1} and \eqref{est-un-lp-2}, changing the power 2 for $p>2$. Moreover we can also prove that $\{u^{(n+1)}(t,x),n\geq0\}$ converges uniformly in $L^p$, denoting this limit by $u(t,x)$. We can check that $u(t,x)$ satisfies Eq. \eqref{mildsolution}. Then it is adapted and satisfies
$$
\sup_{(t,x)\in\mathbb{R}_+\times\mathbb{R}^d}\mathbf{E}(|u(t,x)|^p)<\infty.
$$
The uniqueness can be accomplished by a similar argument.

The key to the continuity is to show that these Picard iterations are mean-square continuous. Then it can be easily extended to $u(t,x)$. In order to show the ideas of the mean-square continuity, we give some steps of the proof for $\{u^{(n+1)}(t,x),n\geq0\}$. As for the time increments, we have, for any $(t,x)\in\mathbb{R}_+\times\mathbb{R}^d$ and $\delta>0$ such that $t+\delta\in\mathbb{R}_+$,
\begin{equation}\label{est-second-moment-1}
\begin{split}
\mathbf{E}&\left[|u^{(n+1)}(t+\delta,x)-u^{(n+1)}(t,x)|^2\right]\\
&\leq \lambda^2\mathbf{E}\left[\left|\int_0^t\int_{\mathbb{R}^d}[G_{t+\delta-u}(x-y)-G_{t-u}(x-y)]
\sigma(u^{(n)}(u,y))W(ds,dy)\right|^2\right]\\
&+\lambda^2\mathbf{E}\left[\left|\int_t^{t+\delta}\int_{\mathbb{R}^d}G_{t+\delta-u}(x-y)
\sigma(u^{(n)}(u,y))W(ds,dy)\right|^2\right]\\
\end{split}
\end{equation}
Using the conditions imposed  on $\sigma$ and \eqref{est-un-lp-1}, we can bound the first term in \eqref{est-second-moment-1} by
$$
C\int_0^{t}du\int_{\mathbb{R}^d}\mu(d\xi)
\left|\mathcal{F}G_{t+\delta-u}(\cdot)(\xi)-\mathcal{F}G_{t-u}(\cdot)(\xi)\right|^2,
$$
which converges to zero as $\delta\downarrow0$ according to \eqref{est-time-1}. The second term in \eqref{est-second-moment-1} can be proved by using the similar arguments by using \eqref{est-time-2}. This proves the right continuity. The left continuity can be proved in the same way.

Concerning the spatial increment, we have, for any $(t,x),(t,z)\in\mathbb{R}_+\times\mathbb{R}^d$,
\begin{equation}\label{est-second-moment-2}
\begin{split}
\mathbf{E}&\left[|u^{(n+1)}(t,x)-u^{(n+1)}(t,z)|^2\right]\\
&\leq C\lambda^2\int_0^{t}du\int_{\mathbb{R}^d}\mu(d\xi)
\left|\mathcal{F}G_{t-u}(x-\cdot)(\xi)-\mathcal{F}G_{t-u}(z-\cdot)(\xi)\right|^2\\
&\leq C\lambda^2\int_0^{t}du\int_{\mathbb{R}^d}\mu(d\xi)|e^{{\rm i}\langle x-z,\xi\rangle}-1|^2
|\mathcal{F}G_{t-u}(z-\cdot)(\xi)|^2.
\end{split}
\end{equation}
Then, thanks to \eqref{est-space-1}, we can prove that the right hand of \eqref{est-second-moment-2} converges to zero as $|x-z|\downarrow0$.

\end{proof}

\qed

\begin{remark}

Let us recall that the SPDE \eqref{sec1-eq1.1} with $\beta=1$ (fractional in space stochastic kinetic equation with factorization of the Laplacian)  has been studied by M\'{a}rquez-Carreras \cite{Carreras2009}.  In this case the Mittag-Leffler function is $E_1(-x)=e^{-x},x\geq0$.

When   $\gamma=0$ and spatial kernel $f(\cdot)$ is the Riesz kernel, then the Eq. \eqref{sec1-eq1.1} reduces to the SPDE studied in Mijena and Nane \cite{MijenaSPA2015}, \cite{MijenaPA2015}. In this reference the authors studied the existence, uniqueness and intermittence of the mild solution for the space-time fractional stochastic partial differential equations \eqref{sec1-eq1.1}.

 For $\gamma=0$ and $\alpha=2$, the SPDE \eqref{sec1-eq1.1} reduces to the classical stochastic heat equation studied by many authors, see e.g. Dalang \cite{DalangEJP1999}

\end{remark}

Now let us make the following assumption on the spectral measure $\mu$ in order to get a precise estimate for the upper bound of the second moment of the mild solution of \eqref{sec1-eq1.1}.
\begin{assumption}\label{assumption-2}
We assume that the spectral measure $\mu$ satisfies
\begin{equation}\label{asymp-nu-1}
\mu(d\xi)\asymp|\xi|^{-\delta}d\xi, \quad {\rm with}\quad 0<\delta<d.
\end{equation}
The symbol $``\asymp"$ means that for every non-negative function $h$ such that the integral in
\eqref{asymp-nu-2} are finite, there exist two positive and finite constants $C$ and $C'$ which may depend on $h$ such that
\begin{equation}\label{asymp-nu-2}
C'\int_{\mathbb{R}^d}h(\xi)|\xi|^{-\delta}d\xi\leq \int_{\mathbb{R}^d}h(\xi)\mu(d\xi)\leq
C\int_{\mathbb{R}^d}h(\xi)|\xi|^{-\delta}d\xi.
\end{equation}
\end{assumption}
\begin{remark}
The Riesz kernel of order $\delta\in(0,d)$ given in Example \ref{example-riesz} obviously satisfies \eqref{asymp-nu-1}. The Bessel kernel given in Example \ref{example-bessel} satisfies \eqref{asymp-nu-1} and the constants in \eqref{asymp-nu-2} are $C=1$ and $C'>0$ depending on $\delta$ and $d$ (see  \cite{TudorXiao2017}).
\end{remark}

We have the following results concerning the upper bound on the second moment of the solution.

\begin{theorem}\label{est-second-moment}
Suppose $0<d-\delta<(\alpha+\gamma)$ and $0<\beta<1$, if the spectral measure $\mu$ associated with the noise $\dot{W}$ satisfies Assumption \ref{assumption-2},
then under the Assumption \ref{assumption-1}, there exist two positive and finite constants $c$ and $c'$ such that
\begin{equation}\label{est-second-moment-pricese}
\sup_{x\in\mathbb{R}^d}\mathbf{E}\left(|u(t,x)|^2\right)\leq c\exp\left\{c'\lambda^{\frac{2(\alpha+\gamma)}{(\alpha+\gamma)-\beta(d-\delta)}}t\right\},
\end{equation}
for all $t>0$.

\end{theorem}

\begin{remark}
This theorem implies that, under some conditions, there exists some positive constant $C$ such that
$$
\limsup_{t\rightarrow\infty}\frac{1}{t}\log \mathbf{E}|u(t,x)|^2\leq C\lambda^{\frac{2(\alpha+\gamma)}{(\alpha+\gamma)-\beta(d-\delta)}},
$$
for any fixed $x\in\mathbb{R}^d$.
\end{remark}

Before giving the proof of Theorem \ref{est-second-moment}, we state  an important lemma needed in the proof of this theorem.

\begin{lemma}\label{cov-3}
Suppose $0<d-\delta<(\alpha+\gamma)$ and $0<\beta<1$, then there exists a positive constant $C$ such that for all $x,y\in\mathbb{R}^d$, we have
$$
\int_{\mathbb{R}^d}\int_{\mathbb{R}^d}G_t(x-z_1)G_t(y-z_2)f(z_1-z_2)dz_1dz_2\leq Ct^{-\frac{\beta(d-\delta)}{\alpha+\gamma}}.
$$
\end{lemma}

\begin{proof}
If we fix $t\in\mathbb{R}_+$, for any $x,y\in\mathbb{R}^d$, then by using \eqref{cov-1}, we have
$$
\int_{\mathbb{R}^d}\int_{\mathbb{R}^d}G_t(x-z_1)G_t(y-z_2)f(z_1-z_2)dz_1dz_2
=\int_{\mathbb{R}^d}\mathcal{F}G_t(x-\cdot)(\xi)
\overline{\mathcal{F}G_t(y-\cdot)(\xi)}\mu(d\xi).
$$
Recall that the spectral measure $\mu$ satisfies \eqref{asymp-nu-1} (i.e. \eqref{asymp-nu-2}) in Assumption \ref{assumption-2}. So according to \eqref{Fourier-Greenfunction}, we have
\begin{equation}
\int_{\mathbb{R}^d}\int_{\mathbb{R}^d}G_t(x-z_1)G_t(y-z_2)f(z_1-z_2)dz_1dz_2
\leq C_{\delta,d}\int_{\mathbb{R}^d}E^2_\beta(-\nu t^\beta|\xi|^\alpha(1+|\xi|^2)^{\gamma/2})|\xi|^{-\delta}d\xi.
\end{equation}
Then by the similar arguments in the proof of Lemma \ref{est-G-square}, based on the estimate on the Mittag-Leffler function (\ref{estimates-E-beta}), we can conclude the proof. \end{proof}

\qed

Now we are ready for giving the proof of Theorem \ref{est-second-moment}. The idea used here is essentially due to \cite{Foondun2015}.

\begin{proof}[Proof of Theorem \ref{est-second-moment}]
Recall the iterated sequences $\{u^{(n)}(t,x),n\geq0, (t,x)\in[0,T]\times\mathbb{R}^d\}$ given by \eqref{iterated-sequences}. Define
\begin{align*}
D_n(t,x)&:=\mathbf{E}\left|u^{(n+1)}(t,x)-u^{(n)}(t,x)\right|^2,\\
H_n(t)&=\sup_{x\in\mathbb{R}^d}D_n(t,x),\\
\Xi(t,y,n)&=|\sigma(u^{(n)}(t,y))-\sigma(u^{(n-1)}(t,y))|.\\
\end{align*}
We will prove the result for $t\in[0,T]$ where $T>0$ is some fixed number.
We now use this notation together with
the covariance formula \eqref{cov-2} and the Assumption \ref{assumption-1} on $\sigma$ to write
$$
D_n(t,x)=\lambda^2\int_0^t\int_{\mathbb{R}^d}\int_{\mathbb{R}^d}G_{t-s}(x-y)G_{t-s}(x-z)
\mathbf{E}(\Xi(s,y,n)\Xi(s,z,n))f(y-z)dydzds.
$$
Now we estimate the expectation on the right hand side using Cauchy-Schwartz
inequality.

\begin{align*}
\mathbf{E} (\Xi(s,y,n)\Xi(s,z,n))&\leq L_\sigma \mathbf{E}\left(|u^{(n)}(s,y)-u^{(n-1)}(s,y)||u^{(n)}(s,z)-u^{(n-1)}(s,z)|\right)\\
&\leq L_\sigma^2\left(\mathbf{E}|u^{(n)}(s,y)-u^{(n-1)}(s,y)|^2\right)^{\frac12}
\left(\mathbf{E}|u^{(n)}(s,z)-u^{(n-1)}(s,z)|^2\right)^{\frac12}\\
&\leq L_\sigma^2\left(D_{n-1}(s,y)D_{n-1}(s,z)\right)^{\frac12}\\
&\leq L_\sigma^2H_{n-1}(s).
\end{align*}

Hence we have for $0<d-\delta<\alpha+\gamma$ by using Lemma \ref{cov-3}
\begin{align*}
D_n(t,x)&\leq\lambda^2L_\sigma^2\int_0^tH_{n-1}(s)
\int_{\mathbb{R}^d}\int_{\mathbb{R}^d}G_{t-s}(x-y)G_{t-s}(x-z)
f(y-z)dydzds\\
&\leq C\lambda^2L_\sigma^2\int_0^tH_{n-1}(s)(t-s)^{-\frac{\beta(d-\delta)}{\alpha+\gamma}}ds.
\end{align*}
We therefore have
$$
H_n(t)\leq C\lambda^2L_\sigma^2\int_0^tH_{n-1}(s)(t-s)^{-\frac{\beta(d-\delta)}{\alpha+\gamma}}ds.
$$
We now note that the integral appearing on the right hand side of the above
inequality is finite when $d-\delta<\frac{\alpha+\gamma}{\beta}$. Hence, by Lemma 3.3 in Walsh \cite{Walsh1986}, the series $\sum_{n=0}^\infty H_n^{\frac12}(t)$
converges uniformly on $[0,T]$. Therefore, the sequence $u^{(n)}(t,x)$ converges
in $L^2$ and uniformly on  $[0,T]\times\mathbb{R}^d$ and the limit satisfies \eqref{mildsolution}. We can
prove uniqueness in a similar way.

We now turn to the proof of the exponential bound. Set
$$
A(t):=\sup_{x\in\mathbb{R}^d}\mathbf{E}|u(t,x)|^2.
$$
We claim that there exist constants $c$ and $c'$ such that for all $t>0$, we have
$$
A(t)\leq c+c'\lambda^2L_\sigma^2\int_0^tA(s)(t-s)^{-\frac{\beta(d-\delta)}{\alpha+\gamma}}ds.
$$
The renewal inequality in Proposition 2.5 in Foondun, Liu and Omaba \cite{FoondunAOP2016} with $\rho=1-\frac{\beta(d-\delta)}{\alpha+\gamma}$, then proves
the exponential upper bound. To prove this claim, we start with the mild
formulation given by \eqref{mildsolution}, then take the second moment to obtain the following
\begin{equation}
\begin{split}
\mathbf{E}&|u(t,x)|^2\\
&=|(\mathcal{G}u_0)_t(x)|^2\\
&+\lambda^2\int_0^t\int_{\mathbb{R}^d}\int_{\mathbb{R}^d}G_{t-s}(x-y)G_{t-s}(x-z)
f(y-z)\mathbf{E}(\sigma(u(s,y))\sigma(u(s,z)))dydzds\\
&:=I_1+I_2.
\end{split}
\end{equation}
Since $u_0$ is bounded, we have $I_1\leq c$. Next we use the Assumption \ref{assumption-1} on $\sigma$ together with H\"{o}lder's inequality to see that
\begin{equation}
\begin{split}
\mathbf{E}(\sigma(u(s,y))\sigma(u(s,z)))&\leq L_\sigma^2 \mathbf{E}(u(s,y)u(s,z))\\
&\leq L_\sigma^2 [\mathbf{E}|u(s,y)|^2]^{\frac12}[\mathbf{E}|u(s,z)|^2]^{\frac12}\\
&\leq L_\sigma^2\sup_{x\in\mathbb{R}^d}\mathbf{E}(|u(s,x)|^2).
\end{split}
\end{equation}
Therefore, using Lemma \ref{cov-3}, the second term $I_2$ is thus bounded  as follows
$$
I_2\leq c\lambda^2L_\sigma^2\int_0^tA(s)(t-s)^{-\frac{\beta(d-\delta)}{\alpha+\gamma}}ds.
$$
Combining the above estimates, we obtain the desired result.

\end{proof}

\qed

Next, we analyze  the H\"{o}lder regularity  of the solution with respect to time and space variables. The next Theorem \ref{thm-holder} extends and improves similar results known for (fractional) stochastic heat equation (e.g. Mijena and Nane \cite{MijenaSPA2015} with $\gamma=0$ in Eq. \eqref{sec1-eq1.1}, Chen and Dalang \cite{ChenSPDE-AC2015}, corresponding to the case $0<\alpha\leq2,\gamma=0$ and $\beta=1$, M\'{a}rquez-Carreras \cite{Carreras2009} with $\beta=1$ in Eq. \eqref{sec1-eq1.1}), and also extends some results for~\eqref{sec1-eq1.1} with Gaussian white noise (e.g. Dalang \cite{DalangEJP1999}). We use a direct method to prove our regularity results in which the Fourier transform and the representation of the Green function (i.e. \eqref{Fourier-Greenfunction} and \eqref{Greenfunction-1}) play a crucial role. We state the result as follows.

\begin{theorem}\label{thm-holder}
Under the Assumption \ref{assumption-1}, assuming that the spectral measure $\mu$ satisfies Hypothesis \ref{hypothesis-2}, then, for every $t,s\in[0,T],T>0$, $x,y\in\mathbb{R}^d$, $p\geq2$ the solution $u(t,x)$ to Eq. \eqref{sec1-eq1.1} satisfies
\begin{equation}
\mathbf{E}\left(|u(t,x)-u(s,y)|^p\right)\leq C\left(|t-s|^{p\chi_1}+|x-y|^{p\chi_2}\right),
\end{equation}
with $0<\chi_1<\min(\beta,\frac12-\beta)$ and $0<\chi_2<\alpha+\gamma-\eta$ if $0<\beta<\frac12$, $0<\chi_1<\frac14$ and $0<\chi_2<\frac{\alpha+\gamma}{2}-\eta$ if $\beta=\frac12$,
and $0<\chi_1<\beta-\frac12$ and $0<\chi_2<\frac{\alpha+\gamma}{2\beta}-\eta$ if $\frac12<\beta<1$.

In particular, the random field $u$ is  $(\chi_1,\chi_2)$-H\"{o}lder continuous with respect to the time and space variables.
\end{theorem}

\begin{proof}
Since the function $(\mathcal{G}u_0)_t(x)=\int_{\mathbb{R}^d}G_t(x-y)u_0(y)dy$ is smooth for any $t>0$, then by Proposition \ref{est-temporal and spatial increments}, \eqref{est-second-moment-1} and \eqref{est-second-moment-2}, we see that for every $p\geq2$ and any $0<T<\infty$, there exists a finite constant $A_{p,T}$ such that
\begin{equation}\label{holder-increments}
\begin{split}
\mathbf{E}\left(|u^{(n)}(t,x)-u^{(n)}(s,y)|^p\right)\leq\left\{
\begin{aligned}
&A_{p,T}\left(|t-s|^{\min(2\beta,1-2\beta)\frac{p}{2}}+|x-y|^{p\chi_2}\right),\quad {\rm if}\quad 0<\beta<\frac12,\\
&A_{p,T}\left(|t-s|^{\frac{p}{4}}+|x-y|^{p\chi_2}\right),\quad {\rm if}\quad \beta=\frac12\\
&A_{p,T}\left(|t-s|^{(2\beta-1)\frac{p}{2}}+|x-y|^{p\chi_2}\right),\quad {\rm if}\quad \frac12<\beta<1.\\
\end{aligned}
\right.
\end{split}
\end{equation}
with $\chi_2\in(0,\alpha+\gamma-\eta)$ if $0<\beta<\frac12$, $\chi_2\in(0,\frac{\alpha+\gamma}{2}-\eta)$ if $\beta=\frac12$ and $\chi_2\in(0,\frac{\alpha+\gamma}{2\beta}-\eta)$ if $\frac12<\beta<1$ simultaneously for all $t,s\in[0,T]$ and $x,y\in\mathbb{R}^d$. The right hand side of this inequality does not depend on $n$. Hence using Fatou's lemma, as $n$ tends to infinity, we get the similar estimates for $u$ which also satisfies \eqref{holder-increments}.
Then the conclusion of Theorem \ref{thm-holder} is a consequence of the Kolmogorov continuity criterion for stochastic processes.
\end{proof}\qed

Let us also make some discussion about the above regularity results.

\begin{remark}
For $\beta$ close to 1, the order of H\"{o}lder regularity of $u(t,x)$ in space is $(\alpha+\gamma)$-times the order of H\"{o}lder continuity to in time. This is coherent with the case in M\'{a}rquez-Carreras \cite{Carreras2009}. When $\gamma=0$ and $\alpha \in(0,2]$, this happens always in the case of the solution of the (fractional) heat equation (with
white noise), see Walsh \cite{Walsh1986}.
\end{remark}
\begin{remark}
If $d=1$, $\alpha=2$ and $\gamma=0$ (so, somehow, the operator $(I-\Delta)^{\frac{\gamma}{2}}(-\Delta)^{\frac\alpha 2}$ reduces to the Laplacian operator $\Delta$ and Eq.\eqref{sec1-eq1.1}
is the standard heat equation), moreover we assume that $\eta$ is close to one-half and $\beta$ is close to 1, we obtain the well-known
regularity of the solution to the heat equation with time-space white noise (which is
H\"{o}lder continuous of order $\frac14$ in time and of order $\frac12$ in space).
\end{remark}

\section{The linear additive noise }\label{sec4}

In the last part  of this work, we focus on the solution of \eqref{mildsolution} with $u_0(x)=0$ and $\sigma=1$. This is the additive noise case and in this situation the solution is Gaussian.  We will study the stationarity of the solution, both ib time and in space. The solution is stationary in space, but not in time.

In the following, we consider
\begin{equation}\label{U}
U(t,x)=\int_0^t\int_{\mathbb{R}^d}G_{t-s}(x-y)W(ds,dy),
\end{equation}
which is the mild solution of \eqref{sec1-eq1.1} when the initial condition $u_0(x)=0, x\in\mathbb{R}^d$ and $\sigma(x)\equiv1,x\in\mathbb{R}$.

\begin{remark}
\begin{enumerate}
\item As mentioned in the introduction, Anh and Leonenko \cite{AnhJSP2001} showed that the presence of the Bessel operator
$-(I-\Delta)^{\gamma/2}$ with $\gamma\geq0$ is essential to have an (asymptotically) stationary solution of SPDE \eqref{sec1-eq1.1}. In fact, the linnear  case requires the condition  $0<\alpha<d/2$ and $\alpha+\gamma>d/2$, that is to say the parameter $\gamma>0$ is necessary.
\item On the other hand, the parameter $\gamma>0$ of the Bessel operator is also useful in determining suitable conditions for the spectral density of the solutions of fractional kinetic equations belonging to $L^1(\mathbb{R}^d)$.
\end{enumerate}
\end{remark}

\begin{theorem}\label{index-x}
Under Hypothesis \ref{hypothesis-2} on the spectral measure $\mu$ associated with $\dot{W}$, then for fixed $t\in\mathbb{R}_+$, the spatial covariance function of $U(t,x)$ given by \eqref{U} is
$$
R_t(x-z)=\int_{\mathbb{R}^d}\mu(d\xi)e^{{\rm i}\langle x-z,\xi\rangle}\int_0^t ds E_\beta^2(-\nu (t-s)^\beta|\xi|^\alpha(1+|\xi|^2)^{\gamma/2}).
$$
In particular, for every $t\in [0,T]$, the process $(U(t,x),  \in \mathbb{R} ^{d})$ is stationary.

\end{theorem}

\begin{proof}
We first calculate the spatial covariance for a fixed $t\in\mathbb{R}_+$. By means of the definition of Fourier transform, change of variable and Fubini's theorem, we obtain, for any $x,z\in\mathbb{R}^d$,
\begin{align*}
\mathbf{E}(U(t,x)\overline{U(t,z)})&=\int_0^t\int_{\mathbb{R}^d}\mathcal{F}G_{t-s}(x-\cdot)(\xi)
\overline{\mathcal{F}G_{t-s}(z-\cdot)(\xi)}\mu(d\xi)ds\\
&=\int_0^t\int_{\mathbb{R}^d}e^{-{\rm i}\langle x-z,\xi\rangle}|\mathcal{F}G_{t-s}(\cdot)(\xi)|^2\mu(d\xi)ds\\
&=\int_{\mathbb{R}^d}e^{-{\rm i}\langle x-z,\xi\rangle}\int_0^tE_\beta^2(-\nu (t-s)^\beta|\xi|^\alpha(1+|\xi|^2)^{\gamma/2})ds\mu(d\xi)\\
&=R_t(x-z),
\end{align*}
where we have used the fact that
\begin{align*}
\int_0^tE_\beta^2&(-\nu (t-s)^\beta|\xi|^\alpha(1+|\xi|^2)^{\gamma/2})ds\\
&\leq \int_0^t\frac{1}{(1+\Gamma(1+\beta)^{-1}(t-s)^\beta|\xi|^\alpha(1+|\xi|^2)^{\gamma/2})^2}ds\\
&\leq t.
\end{align*}
Hence, for fixed $t\in\mathbb{R}_+$, the process $U(t,x),x\in\mathbb{R}^d$ is a Gaussian field that has zero mean, stationary increments, and a continuous covariance function. \end{proof}\qed

\begin{remark}
From the above result, one can obtain the spectral density of the process $x\to U(t,x)$. Indeeed, its spectral density $f_{t}(\lambda)$ is given by
$$f_{t}(\lambda)= \int_0^t ds E_\beta^2(-\nu (t-s)^\beta|\xi|^\alpha(1+|\xi|^2)^{\gamma/2}) g(\xi)$$
where $g$ is the density of $\mu$ with respect to the Lebesque masure.

\end{remark}

The next result shows that the process \eqref{U} is not stationary in time but, when $t$ tends to infinity, it converges to a stationary process.

\begin{theorem}\label{index-t}
Under Hypothesis \ref{hypothesis-2} on the spectral measure $\mu$ associated with $\dot{W}$, assuming $1/2<\beta<1$, then for $t\in\mathbb{R}_+$, $\tau\in\mathbb{R}$ such that $t+\tau\in\mathbb{R}_+$, $x,z\in\mathbb{R}^d$, the asymptotic homogeneous spatial-temporal covariance function of $U(t+\tau,x)$ and $U(t,z)$ is
\begin{align*}
R(\tau,x-z)=\int_{\mathbb{R}^d}&
\frac{e^{{\rm i}\langle x-z,\xi\rangle}}{\beta(\nu|\xi|^\alpha(1+|\xi|^2)^{\gamma/2})^{1/\beta}}
\int_{\nu\tau^\beta|\xi|^\alpha(1+|\xi|^2)^{\gamma/2}}^{\infty}x^{1/\beta-1}E_\beta(-x)\\
&\hspace{1cm}\cdot E_\beta\left(-\nu\left(\left(\frac{x}{\nu|\xi|^\alpha(1+|\xi|^2)^{\gamma/2}}\right)
^{1/\beta}-\tau\right)^\beta|\xi|^\alpha(1+|\xi|^2)^{\gamma/2}\right)dx\mu(d\xi).
\end{align*}
 Moreover $U(\cdot,x)$ is asymptotically in time an index-$\left(\beta-\frac12\right)$ Gaussian field.

\end{theorem}

\begin{proof}
For $t,\tau\in\mathbb{R}_+$ (for $\tau\in\mathbb{R}_{-}$ such that $t+\tau\in\mathbb{R}_+$, we argue similarly), and $x,z\in\mathbb{R}^d$, we have
\begin{align*}
&\mathbf{E}\left(U(t+\tau,x)\overline{U(t,z)}\right)\\
&=\int_0^tds\int_{\mathbb{R}^d}\mu(d\xi)
\mathcal{F}G_{t+\tau-s}(x-\cdot)(\xi)\overline{\mathcal{F}G_{t-s}(z-\cdot)(\xi)}\\
&=\int_0^tds\int_{\mathbb{R}^d}\mu(d\xi)e^{{\rm i}\langle x-z,\xi\rangle}
\mathcal{F}G_{t+\tau-s}(x-\cdot)(\xi)\mathcal{F}G_{t-s}(z-\cdot)(\xi)\\
&=\int_0^tds\int_{\mathbb{R}^d}\mu(d\xi)e^{{\rm i}\langle x-z,\xi\rangle}
E_\beta(-\nu(t+\tau-s)^\beta|\xi|^\alpha(1+|\xi|^2)^{\gamma/2})
E_\beta(-\nu(t-s)^\beta|\xi|^\alpha(1+|\xi|^2)^{\gamma/2})\\
&=\int_{\mathbb{R}^d}e^{{\rm i}\langle x-z,\xi\rangle}
\int_0^t E_\beta(-\nu(t+\tau-s)^\beta|\xi|^\alpha(1+|\xi|^2)^{\gamma/2})
E_\beta(-\nu(t-s)^\beta|\xi|^\alpha(1+|\xi|^2)^{\gamma/2})ds\mu(d\xi)\\
\end{align*}
Next let us calculate the above integral with respect to $s$. In fact, with the change of variable $x=\nu(t+\tau-s)^\beta|\xi|^\alpha(1+|\xi|^2)^{\gamma/2}$ we have
\begin{align*}
\int_0^t& E_\beta(-\nu(t+\tau-s)^\beta|\xi|^\alpha(1+|\xi|^2)^{\gamma/2})
E_\beta(-\nu(t-s)^\beta|\xi|^\alpha(1+|\xi|^2)^{\gamma/2})ds\\
&=\frac{1}{\beta(\nu|\xi|^\alpha(1+|\xi|^2)^{\gamma/2})^{1/\beta}}
\int_{\nu\tau^\beta|\xi|^\alpha(1+|\xi|^2)^{\gamma/2}}^{\nu(t+\tau)
^\beta|\xi|^\alpha(1+|\xi|^2)^{\gamma/2}}E_\beta(-x)x^{1/\beta-1}\\
&\hspace{2cm}\cdot E_\beta\left(-\nu\left(\left(\frac{x}{\nu|\xi|^\alpha(1+|\xi|^2)^{\gamma/2}}\right)
^{1/\beta}-\tau\right)^\beta|\xi|^\alpha(1+|\xi|^2)^{\gamma/2}\right)dx.
\end{align*}
Moreover, as $t\rightarrow\infty$, we obtain
\begin{align*}
&\lim_{t\rightarrow\infty}\int_0^t E_\beta(-\nu(t+\tau-s)^\beta|\xi|^\alpha(1+|\xi|^2)^{\gamma/2})
E_\beta(-\nu(t-s)^\beta|\xi|^\alpha(1+|\xi|^2)^{\gamma/2})ds\\
&=\frac{1}{\beta(\nu|\xi|^\alpha(1+|\xi|^2)^{\gamma/2})^{1/\beta}}
\int_{\nu\tau^\beta|\xi|^\alpha(1+|\xi|^2)^{\gamma/2}}^{\infty}x^{1/\beta-1}E_\beta(-x)\\
&\hspace{2cm}\cdot E_\beta\left(-\nu\left(\left(\frac{x}{\nu|\xi|^\alpha(1+|\xi|^2)^{\gamma/2}}\right)
^{1/\beta}-\tau\right)^\beta|\xi|^\alpha(1+|\xi|^2)^{\gamma/2}\right)dx,
\end{align*}
which is finite. Because we have
\begin{align*}
\int_{\nu\tau^\beta|\xi|^\alpha(1+|\xi|^2)^{\gamma/2}}^{\infty}&x^{1/\beta-1}E_\beta(-x) E_\beta\left(-\nu\left(\left(\frac{x}{\nu|\xi|^\alpha(1+|\xi|^2)^{\gamma/2}}\right)
^{1/\beta}-\tau\right)^\beta|\xi|^\alpha(1+|\xi|^2)^{\gamma/2}\right)dx\\
&\leq C_{\xi,\nu,\beta,\alpha,\gamma,\tau}\int_{\nu\tau^\beta|\xi|^\alpha(1+|\xi|^2)}^{\infty}
x^{1/\beta-3}dx,
\end{align*}
which is finite when $1/2<\beta<1$ (i.e. $\frac1\beta-2<0$).

We now tackle the second part of this theorem. We assume that $x\in\mathbb{R}^d$, $t\in\mathbb{R}_+$ and $\tau\in\mathbb{R}_+$ are small (the negative case is similar). Then, from \eqref{est-time-1} and \eqref{est-time-2}, when $1/2<\beta<1$  we have
\begin{align*}
\mathbf{E}&\left(|U(t+\tau,x)-U(t,x)|^2\right)\\
&=\int_0^{t}ds\int_{\mathbb{R}^d}\mu(d\xi)
\left|\mathcal{F}G_{t+\tau-s}(x-\cdot)(\xi)-\mathcal{F}G_{t-s}(x-\cdot)(\xi)\right|^2\\
&\hspace{1.5cm}+\int_{t}^{t+\tau}ds\int_{\mathbb{R}^d}\mu(d\xi)
\left|\mathcal{F}G_{t+\tau-s}(x-\cdot)(\xi)\right|^2\\
&\leq C\left(|\tau|^{2\beta}+|\tau|^{2\beta-1}\right).
\end{align*}
Then we can complete the proof of the second part.

\end{proof}

\qed

\section*{Appendix}

We will prove Proposition \ref{est-temporal and spatial increments} in this appendix.
\begin{proof}[Proof of Proposition \ref{est-temporal and spatial increments}]
For any $t,t'\in\mathbb{R}_+$ such that $t'<t$ and $x\in\mathbb{R}^d$,
by using the fact
$$
\mathcal{F}G_{t-s}(x-\cdot)(\xi)-\mathcal{F}G_{t'-s}(x-\cdot)(\xi)
=e^{{\rm i}\langle\xi,x\rangle}(\mathcal{F}G_{t-s}(\cdot)(\xi)
-\mathcal{F}G_{t'-s}(\cdot)(\xi)),
$$
then from \eqref{Fourier-Greenfunction} and the absolute convergence of the series in \eqref{Def-E-1}, one gets
\begin{equation}\label{fourier-difference G}
\begin{split}
\mathcal{F}G_{t-s}(\cdot)(\xi)
-\mathcal{F}G_{t'-s}(\cdot)(\xi)&=\sum_{k=0}^\infty
\frac{(-\nu|\xi|^\alpha(1+|\xi|^2)^{\gamma/2})^k}{\Gamma(1+k\beta)}
\left(|t-s|^{k\beta}-|t'-s|^{k\beta}\right)\\
&\leq \sum_{k=0}^\infty
\frac{(-\nu|\xi|^\alpha(1+|\xi|^2)^{\gamma/2})^k}{\Gamma(1+k\beta)}\left( |t-s|^{k}-|t'-s|^{k}\right)^\beta\\
&\leq |t-t'|^\beta\sum_{k=1}^\infty
\frac{k^\beta(-\nu|\xi|^\alpha(1+|\xi|^2)^{\gamma/2})^k}{\Gamma(1+k\beta)}t^{(k-1)\beta},
\end{split}
\end{equation}
for all $t,t'\in \mathbb{R}_+$, where the last inequality follows from the mean value theorem. Furthermore, since the series in \eqref{Def-E-1} is absolutely convergence, then the series in the last inequality in \eqref{fourier-difference G} can be bounded as follows with $0<\beta<1$
\begin{equation}\label{fourier-difference G-1}
\begin{split}
\sum_{k=1}^\infty
\frac{k^\beta(-\nu|\xi|^\alpha(1+|\xi|^2)^{\gamma/2})^k}{\Gamma(1+k\beta)}t^{(k-1)\beta}&
=t^{-\beta}\sum_{k=1}^\infty\frac{k^\beta(-\nu|\xi|^\alpha(1+|\xi|^2)^{\gamma/2} t^\beta)^k}{\Gamma(1+k\beta)}\\
&\leq-\nu  t^{-\beta}|\xi|^\alpha(1+|\xi|^2)^{\gamma/2} \sum_{k=1}^\infty\frac{k(-\nu|\xi|^\alpha(1+|\xi|^2)^{\gamma/2} t^\beta)^{k-1}}{\Gamma(1+k\beta)}\\
&\leq-\nu t^{-\beta}|\xi|^\alpha(1+|\xi|^2)^{\gamma/2} \sum_{k=1}^\infty\frac{(-\nu|\xi|^\alpha(1+|\xi|^2)^{\gamma/2} t^\beta)^{k-1}}{\Gamma(1+k\beta)}\\
&\leq\frac{t^{-\beta}}{1+\Gamma(1+\beta)^{-1}\nu t^\beta|\xi|^\alpha (1+|\xi|^2)^{\gamma/2}}.
\end{split}
\end{equation}
Then we have
\begin{equation}\label{fourier-difference G-2}
\begin{split}
\int_0^{t'}ds&\int_{\mathbb{R}^d}\mu(d\xi)
\left|\mathcal{F}G_{t-s}(x-\cdot)(\xi)
-\mathcal{F}G_{t'-s}(x-\cdot)(\xi)\right|^2\\
&\leq |t-t'|^{2\beta}\int_0^{t'}ds\int_{\mathbb{R}^d}\mu(d\xi)
\frac{t^{-2\beta}}{(1+\Gamma(1+\beta)^{-1}\nu t^\beta|\xi|^\alpha (1+|\xi|^2)^{\gamma/2})^2}\\
&:=A_1+A_2,
\end{split}
\end{equation}
with
$$
A_1=|t-t'|^{2\beta}\int_0^{t'}ds\int_{|\xi|\leq1}\mu(d\xi)
\frac{t^{-2\beta}}{(1+\Gamma(1+\beta)^{-1}\nu t^\beta|\xi|^\alpha (1+|\xi|^2)^{\gamma/2})^2},
$$
and
$$
A_2=|t-t'|^{2\beta}\int_0^{t'}ds\int_{|\xi|>1}\mu(d\xi)
\frac{t^{-2\beta}}{(1+\Gamma(1+\beta)^{-1}\nu t^\beta|\xi|^\alpha (1+|\xi|^2)^{\gamma/2})^2}.
$$
So with \eqref{est-mu}, we have
$$
A_1\leq t^{1-2\beta}\int_{|\xi|\leq1}\mu(d\xi)|t-t'|^{2\beta}.
$$
and
\begin{align*}
A_2&\leq t^{-2\beta}\int_{|\xi|>1}N_{t'}(\xi)\mu(d\xi)|t-t'|^{2\beta}
\end{align*}
Then combining the above estimates for $A_1$, $A_2$ and Proposition \ref{nt-xi-lower and upper}, we can conclude the proof of \eqref{est-time-1}. Next we can follow the similar arguments to prove \eqref{est-time-2} which will be divided into three cases. Firstly with $0<\beta<\frac12$, we have
\begin{align*}
\int_{t'}^t&ds\int_{\mathbb{R}^d}\mu(d\xi)
\left|\mathcal{F}G_{t-s}(x-\cdot)(\xi)\right|^2\\
&=\int_{t'}^tds\int_{\mathbb{R}^d}\mu(d\xi)|e^{{\rm i}\langle\xi,x\rangle}|^2\left|E_\beta(-\nu(t-s)^\beta|\xi|^\alpha (1+|\xi|^2)^{\gamma/2})\right|^2\\
&\leq \int_{\mathbb{R}^d}\mu(d\xi)\int_{t'}^tds\left(\frac{1}{1+\Gamma(1+\beta)^{-1}
\nu(t-s)^\beta|\xi|^\alpha (1+|\xi|^2)^{\gamma/2}}\right)^2\\
&\leq \int_{|\xi|\leq1}\mu(d\xi)|t-t'|\\
&\hspace{1cm}+\frac{\Gamma(1+\beta)^22^{\alpha+\gamma}}{\nu^2(1-2\beta)}
\int_{|\xi|>1}\left(\frac{1}{1+|\xi|^2}\right)
^{\alpha+\gamma}\mu(d\xi)|t-t'|^{1-2\beta}\\
&=C_{3.2}|t-t'|^{1-2\beta}.
\end{align*}
If $\frac12<\beta<1$, one obtains with Fubini's theorem
\begin{align*}
\int_{t'}^tds&\int_{\mathbb{R}^d}\mu(d\xi)
\left|\mathcal{F}G_{t-s}(x-\cdot)(\xi)\right|^2\\
&=\int_{\mathbb{R}^d}\mu(d\xi)\int_{t'}^t
\left|E_\beta(-\nu(t-s)^\beta|\xi|^\alpha (1+|\xi|^2)^{\gamma/2})\right|^2ds\\
&\leq \int_{\mathbb{R}^d}\mu(d\xi)\int_{t'}^t\left(\frac{1}{1+\Gamma(1+\beta)^{-1}
\nu(t-s)^\beta|\xi|^\alpha (1+|\xi|^2)^{\gamma/2}}\right)^2ds.
\end{align*}
Denote by
\begin{align*}
M_{t,t'}(\xi)&:=\int_{t'}^t\left(\frac{1}{1+\Gamma(1+\beta)^{-1}
\nu(t-s)^\beta|\xi|^\alpha (1+|\xi|^2)^{\gamma/2}}\right)^2ds\\
&=\int_{0}^{t-t'}\left(\frac{1}{1+\Gamma(1+\beta)^{-1}
\nu u^\beta|\xi|^\alpha (1+|\xi|^2)^{\gamma/2}}\right)^2du
\end{align*}
Then with the change of variable $x=\Gamma(1+\beta)^{-1}
\nu u^\beta|\xi|^\alpha (1+|\xi|^2)^{\gamma/2}$, we have
\begin{align*}
M_{t,t'}(\xi)=\frac{1}{\beta}\left(\frac{\Gamma(1+\beta)}{\nu|\xi|^\alpha (1+|\xi|^2)^{\gamma/2}}\right)^{\frac1\beta}\int_{0}^{\Gamma(1+\beta)^{-1}
\nu (t-t')^\beta|\xi|^\alpha (1+|\xi|^2)^{\gamma/2}}x^{\frac1\beta-1}\left(1+x\right)^{-2}dx.
\end{align*}
For $|\xi|\leq1$, we have the following
\begin{equation}\label{est-M-1}
\begin{split}
M_{t,t'}(\xi)1_{\{|\xi|\leq1\}}&\leq\frac{1}{\beta}\left(\frac{\Gamma(1+\beta)}{\nu|\xi|^\alpha (1+|\xi|^2)^{\gamma/2}}\right)^{\frac1\beta}\int_{0}
^{\frac{\nu|\xi|^\alpha (1+|\xi|^2)^{\gamma/2}}{\Gamma(1+\beta)}(t-t')^\beta
}x^{\frac1\beta-1}dx\\
&=|t-t'|.\\
\end{split}
\end{equation}
Now we will estimate $M_{t,t'}(\xi)$ when $|\xi|>1$. In fact, with $0<\eta<\frac{\alpha+\gamma}{2\beta}$, we have
\begin{equation}\label{est-M-2}
\begin{split}
M_{t,t'}(\xi)1_{\{|\xi|>1\}}
&\leq\frac{1}{\beta}\left(\frac{\Gamma(1+\beta)}{\nu|\xi|^\alpha (1+|\xi|^2)^{\gamma/2}}\right)^{\frac1\beta}\int_{0}
^{\frac{\nu|\xi|^\alpha (1+|\xi|^2)^{\gamma/2}}{\Gamma(1+\beta)}(t-t')^\beta
}(1+x)^{\frac1\beta-3}dx\\
&\leq\frac{1}{2\beta-1}\left(\frac{\Gamma(1+\beta)}
{\nu|\xi|^\alpha (1+|\xi|^2)^{\gamma/2}}\right)^{\frac1\beta}\\
&\leq \frac{1}{2\beta-1}\left(\frac{\Gamma(1+\beta)}{\nu}\right)
^{\frac1\beta}2^{\frac{\alpha+\gamma}{2\beta}}
\left(\frac{1}{1+|\xi|^2}\right)^{\frac{\alpha+\gamma}{2\beta}}.
\end{split}
\end{equation}
Then combining the estimates \eqref{est-M-1} and \eqref{est-M-2}, we can obtain \eqref{est-time-2}.

If $\beta=\frac12$, by using the similar argument above, we have
\begin{align*}
\int_{t'}^tds&\int_{\mathbb{R}^d}\mu(d\xi)
\left|\mathcal{F}G_{t-s}(x-\cdot)(\xi)\right|^2\\
&=\int_{\mathbb{R}^d}\mu(d\xi)\int_{t'}^t
\left|E_\beta(-\nu(t-s)^\beta|\xi|^\alpha (1+|\xi|^2)^{\gamma/2})\right|^2ds\\
&\leq \int_{|\xi|\leq1}\mu(d\xi)|t-t'|\\
&+2\left(\frac{\Gamma(3/2)}{\nu|\xi|^\alpha (1+|\xi|^2)^{\gamma/2}}\right)^2\int_{|\xi|>1}\mu(d\xi)\int_{0}^{\Gamma(3/2)^{-1}\nu|t-t'|^{\frac12}
|\xi|^\alpha (1+|\xi|^2)^{\gamma/2}}x(1+x)^{-2}dx.
\end{align*}
and we have
\begin{align*}
&2\left(\frac{\Gamma(3/2)}{\nu|\xi|^\alpha (1+|\xi|^2)^{\gamma/2}}\right)^2\int_{0}^{\Gamma(3/2)^{-1}\nu|t-t'|^{\frac12}
|\xi|^\alpha (1+|\xi|^2)^{\gamma/2}}x(1+x)^{-2}dx\\
&\leq2\left(\frac{\Gamma(3/2)}{\nu|\xi|^\alpha (1+|\xi|^2)^{\gamma/2}}\right)^2\int_{0}^{\Gamma(3/2)^{-1}\nu|t-t'|^{\frac12}
|\xi|^\alpha (1+|\xi|^2)^{\gamma/2}}(1+x)^{-1}dx\\
&=2\left(\frac{\Gamma(3/2)}{\nu|\xi|^\alpha (1+|\xi|^2)^{\gamma/2}}\right)^2\ln\left(1+\Gamma(3/2)^{-1}\nu|t-t'|^{\frac12}
|\xi|^\alpha (1+|\xi|^2)^{\gamma/2}\right)\\
&\leq2^{1+\frac{\alpha+\gamma}{2}}\frac{\Gamma(3/2)}{\nu}\left(\frac{1}{1+|\xi|^2}\right)^{\frac{\alpha+\gamma}{2}}|t-t'|^{\frac12}.
\end{align*}
Then we have that
\begin{align*}
\int_{t'}^tds&\int_{\mathbb{R}^d}\mu(d\xi)
\left|\mathcal{F}G_{t-s}(x-\cdot)(\xi)\right|^2\\
&\leq \left(\int_{|\xi|\leq1}\mu(d\xi)|t-t'|^{\frac12}+2^{1+\frac{\alpha+\gamma}{2}}\frac{\Gamma(3/2)}{\nu}\int_{|\xi|>1}\left(\frac{1}{1+|\xi|^2}\right)^{\frac{\alpha+\gamma}{2}}\mu(d\xi)\right)|t-t'|^{\frac12}\\
\end{align*}

Since
$$
\mathcal{F}G_{t-s}(x-\cdot)(\xi)-\mathcal{F}G_{t-s}(x'-\cdot)(\xi)
=(e^{{\rm i}\langle\xi,x\rangle}-e^{{\rm i}\langle\xi,x'\rangle})\mathcal{F}G_{t-s}(\cdot)(\xi),
$$
we have
\begin{align*}
\int_0^tds&\int_{\mathbb{R}^d}\mu(d\xi)
\left|\mathcal{F}G_{t-s}(x-\cdot)(\xi)-\mathcal{F}G_{t-s}(x'-\cdot)(\xi)\right|^2\\
&=\int_0^tds\int_{\mathbb{R}^d}\mu(d\xi)
\left|e^{{\rm i}\langle\xi,x\rangle}-e^{{\rm i}\langle\xi,x'\rangle}\right|^2\left|\mathcal{F}G_{t-s}(\cdot)(\xi)\right|^2\\
&:=B_1+B_2,
\end{align*}
with
\begin{align*}
B_1&=\int_0^tds\int_{|\xi|\leq1}\mu(d\xi)
\left|e^{{\rm i}\langle\xi,x\rangle}-e^{{\rm i}\langle\xi,x'\rangle}\right|^2\left|\mathcal{F}G_{t-s}(\cdot)(\xi)\right|^2\\
B_2&=\int_0^tds\int_{|\xi|>1}\mu(d\xi)
\left|e^{{\rm i}\langle\xi,x\rangle}-e^{{\rm i}\langle\xi,x'\rangle}\right|^2\left|\mathcal{F}G_{t-s}(\cdot)(\xi)\right|^2.\\
\end{align*}
The first term $B_1$ is easy and can be studied in the same way for any $0<\beta<1$. Indeed, the fact that the Fourier transform of Green function $G$ given by \eqref{Fourier-Greenfunction} is bounded by 1, the mean value theorem, and property \eqref{est-mu} imply that
\begin{equation}\label{est-b1}
\begin{split}
B_1&\leq C\int_0^tds\int_{|\xi|\leq1}\mu(d\xi)|\langle x-x',\xi\rangle|^2\\
&\leq Ct\int_{|\xi|\leq1}\mu(d\xi)|x-x'|^2.
\end{split}
\end{equation}
The other term $B_2$ is a little involved. We distinguish three cases depending on the values of $\beta$. We first study the case $0<\beta<\frac12$. Let $0<\rho_1<\alpha+\gamma-\eta$.  Applying the mean theorem, Fubini's theorem, the fact $1-e^{-x}\leq1$ for all $x>0$ and the Hypothesis \ref{hypothesis-2},
then we have
\begin{equation}\label{est-b2-leq}
\begin{split}
B_2&=\int_0^tds\int_{|\xi|>1}\mu(d\xi)
\left|e^{{\rm i}\langle\xi,x\rangle}-e^{{\rm i}\langle\xi,x'\rangle}\right|^2\left|\mathcal{F}G_{t-s}(\cdot)(\xi)\right|^2\\
&\leq4\int_0^tds\int_{|\xi|>1}\mu(d\xi)\left|\frac12(e^{{\rm i}\langle\xi,x\rangle}-e^{{\rm i}\langle\xi,x'\rangle})\right|^{2\rho_3}\left|\mathcal{F}G_{t-s}(\cdot)(\xi)\right|^2\\
&\leq C\int_0^tds\int_{|\xi|>1}\mu(d\xi)|\xi|^{2\rho_1}|x-x'|^{2\rho_1}
|E_\beta(-\nu(t-s)^\beta|\xi|^\alpha(1+|\xi|^2)^{\gamma/2})|^2\\
&\leq Ct^{1-2\beta} \frac{\Gamma(1+\beta)^{2}}{(1-2\beta)\nu^{2}}2^{\alpha+\gamma-\rho_1}
\int_{|\xi|>1}\left(\frac{1}{1+|\xi|^2}\right)
^{\alpha+\gamma-\rho_1}\mu(d\xi)|x-x'|^{2\rho_1}.\\
\end{split}
\end{equation}
For the critical case $\beta=\frac12$, by choosing $0<\rho_2<\frac{\alpha+\gamma}{2}-\eta$, then we have that
\begin{equation}\label{est-b2-equal}
\begin{split}
B_2&\leq C|x-x'|^{2\rho_2}\int_{|\xi|>1}\mu(d\xi)|\xi|^{2\rho_2}\int_0^t
|E_\beta(-\nu(t-s)^\beta|\xi|^\alpha(1+|\xi|^2)^{\gamma/2})|^2ds\\
&\leq 2Ct^{\frac12}|x-x'|^{2\rho_2}\Gamma(3/2)\nu^{-1}\int_{|\xi|>1}\mu(d\xi)\frac{|\xi|^{2\rho_2}}{|\xi|^\alpha(1+|\xi|^2)^{\gamma/2})}\\
&\leq C 2^{1+\frac{\alpha+\gamma}{2}-\rho_2}t^{\frac12}|x-x'|^{2\rho_2}\Gamma(3/2)\nu^{-1}\int_{|\xi|>1}\left(\frac{1}{1+|\xi|^2}\right)^{\frac{\alpha+\gamma}{2}-\rho_2}\mu(d\xi).
\end{split}
\end{equation}

On the other hand, when $\frac12<\beta<1$, let $0<\rho_3<\frac{\alpha+\gamma}{2\beta}-\eta$, then the similar arguments yield that
\begin{equation}\label{est-b2-geq}
\begin{split}
B_2&\leq C\int_0^tds\int_{|\xi|>1}\mu(d\xi)|\xi|^{2\rho_3}|x-x'|^{2\rho_3}
|E_\beta(-\nu(t-s)^\beta|\xi|^\alpha(1+|\xi|^2)^{\gamma/2})|^2\\
&\leq \frac{C}{(2\beta-1)}\left(\frac{\Gamma(1+\beta)}{\nu}\right)
^{\frac1\beta}
2^{\frac{\alpha+\gamma}{2\beta}-\rho_3}
\int_{|\xi|>1}\left(\frac{1}{1+|\xi|^2}\right)
^{\frac{\alpha+\gamma}{2\beta}-\rho_3}\mu(d\xi)\cdot|x-x'|^{2\rho_3}.\\
\end{split}
\end{equation}

Then we can conclude the proof of \eqref{est-space-1} by combining \eqref{est-b1}, \eqref{est-b2-geq}, \eqref{est-b2-leq} and \eqref{est-b2-equal}.

\end{proof}

\qed

\end{document}